\documentclass[12pt,leqno]{amsart}
\usepackage{amsmath,amssymb,amsthm}
\usepackage[pdftex]{graphicx}
\usepackage{qtree,tcolorbox,tikz}

\begin{document}

\title[Irreducible representations of $B_n$]{Irreducible representations of Braid Group $B_n$ of dimension $n+1$}

\author{Inna Sysoeva}

\email{inna@math.binghamton.edu}
\date{May 26, 2020}

\begin{abstract}
We prove that there are no irreducible representations of $B_n$ of dimension $n+1$ for $n\geqslant 10.$
\end{abstract}

\maketitle

\section{Introduction}

The irreducible complex representations of braid group $B_n$ on $n$ strings of dimension less than or equal to $n$ have been completely classified about 20 years ago. Formanek (see \cite{Formanek}) had classified all the irreducible repesentations of dimension of less than or equal to $(n-1).$  In \cite{S}, the author provided the classification of irreducible representations of dimension $n$ for $n\geqslant 9,$ and the classification for the small values of $n$ was completed in   \cite{Lee}  and  \cite{FLSV}.

  Apart from the number of exceptions for small values of $n,$ all irreducible representations of $B_n$ are either one-dimensional, or have dimensions $(n-2)$, $(n-1),$ or $n.$ All $(n-2)-$ dimensional representations are equivalent to a tensor product of a one-dimensional representation and an $(n-2)-$dimensional composition factor of a specialization of the reduced Burau representation (see \cite{Burau}, \cite{Jones87}, \cite{Formanek}); $(n-1)-$dimensional representations are equivalent to  a tensor product of a one-dimensional representation and a specialization of the reduced Burau representation, and $n$-dimensional  representations are equivalent to  a tensor product of a one-dimensional representation and a specialization of the standard representation
(see \cite{TongYangMa}). 

Another famous representation (whose dimension is $\frac{n(n-1)}{2}$) is  the Lawrence-Krammer-Bigelow representation (see \cite{Lawrence}). Its irreducibility   for generic values of the parameters was first proven by  Zinno in \cite{Zinno} and, as Zinno mentions in his paper,
V. F. R. Jones also discovered the same  result, simultaneously and independently. 

However, the problem of the systematic classification of the irreducible representations, even between  dimensions $n$ and  $\frac{n(n-1)}{2},$ is wide open. To the best of my knowledge, up to this point no such classification was found. No irreducible representations of dimensions between $n$ and $\frac{n(n-1)}{2}$ were found for $n$ large enough, and no results of non-existence for the irreducible representations are known. The most notable step towards such classification was published by Larsen and Rowell in \cite{LarsenRowell}, where, among other results, they  had proven that there are no irreducible {\it unitary} representations of $B_n$ between dimensions $n+1$ and $2n-9$ (inclusive) for $n$ large enough. In particular, there are no irreducible {\it unitary} representations of dimension $n+1$ for $n\geqslant 16.$

The main result of this paper (Theorem 6.1) is that $B_n$ has no irreducible representations of dimension $n+1$ for $n\geqslant 10.$ We do not claim that this lower bound is sharp. As Larsen and Rowell pointed out in their paper, the actual sharp lower bound is at least 8, since $B_7$ has an 8-dimensional unitary representation (see Jones, \cite{Jones}). 

To achieve our  result, we first prove  that for $n\geqslant 10,$ for any irreducible representation of $B_n$ of dimension not exceeding $2n-9,$ the image of each group generator $\rho(\sigma_i)$  must have an eigenvalue of high multiplicity; in other words, the representation is equivalent to the tensor product of a one-dimensional representation and an irreducible representation of small corank (The Reduction Theorem, Theorem 3.8). This theorem  by itself is a separate tool which may be used in the classification of the irreducible representations of dimension up to $2n-9.$

The paper is organized as follows. In section 2 we give the definitions and recall preliminary results. In section 3 the Reduction Theorem is proven. In section 4 we review the definition of the {\it friendship graph} which was introduced in \cite{S}. We will briefly review the results from \cite{S} related to the friendship graphs, and prove some new results that will be used to classify the irreducible  representations of dimension $n+1.$ In section 5 we prove that there are no irreducible representations of dimension $r\geqslant n+1$ and corank 3 for $n\geqslant 10,$ and we  summarize the main results of this paper in section 6, Theorems 6.1 and 6.2.\\

{\bf Acknowledgements.} I would like to take this opportunity to express my deep gratitude to Edward Formanek, who spiked my interest in this subject long time ago. I also would like to thank my husband Alexander Borisov who convinced me to return to this work after significant number of  years  despite challenging life circumstances. This work would never have been completed without his persistence and moral support.\\

\section{Notations and Preliminary Results}

Let $B_n$ be the braid group on $n$ strings. As an abstract group, it has the following presentation:\\

$B_n=\left\langle \begin{tabular}{l|l}
&$\sigma_i \sigma_j=\sigma_j\sigma_i$ for $|i-j|\geqslant 2,$\\
$\sigma_1,\ \sigma_2,  \dots , \sigma_{n-1}$ &\\
& $\sigma_i \sigma_{i+1}\sigma_i=\sigma_{i+1}\sigma_i \sigma_{i+1}$ for $1\leqslant i\leqslant n-2$
\end{tabular}\right\rangle, $\\

where $\sigma_1,\ \sigma_2,  \dots , \sigma_{n-1}$ are called  the standard generators.\\

Since $B_2\cong \mathbb{Z},$ we will assume throughout the paper that $n\geqslant 3,$ unless specified otherwise.\\

In the following lemma we are going to recall the results originally due to W.-L. Chow \cite{Chow} and F.A. Garside \cite{Garside}.\\

{\bf Lemma 2.1.} Let $B_n$ be a braid group on $n$ strings, $n\geqslant 3.$  \\

Let $\tau=\sigma_1\sigma_2\dots \sigma_{n-1},$ and let  $\sigma_{0}=\tau\sigma_{n-1}\tau^{-1}.$\\

Let $\Delta=(\sigma_{n-1}\sigma_{n-2}\dots\sigma_1)(\sigma_{n-1}\sigma_{n-2}\dots \sigma_2)\dots(\sigma_{n-1}\sigma_{n-2})(\sigma_{n-1})$ \\(positive half-twist).\\

Then:\\

(a) (\cite{Chow}, Equation (6)) $\sigma_{i+1}=\tau\sigma_{i}\tau^{-1}$ for $1\leqslant i\leqslant n-2;$

(b)(\cite{Chow}, Theorem III): The center of $B_n,$  $Z(B_n)=<\tau ^n>;$

(c) $\sigma_{1}=\tau\sigma_{0}\tau^{-1};$ 

 (d) (\cite{Garside}, Theorem 7) $Z(B_n)=<\Delta ^2>;$ 

(e) (\cite{Garside}, Lemma 2(ii)) $\sigma_{n-i}=\Delta \sigma_i \Delta^{-1},  \ i=1, \ 2, \dots , n-1.$\\

{\it Proof of }(c). It easily follows from (a), (b) and the definition of $\sigma_0.$

 \begin{flushright}$\square$\end{flushright}

In view of Lemma 2.1, we can add  a redundant generator $\sigma_0$ to the set of standard generators of $B_n$ to get the following presentation:\\
 
$B_n=\left\langle \begin{tabular}{l|l}
&$\sigma_i \sigma_j=\sigma_j\sigma_i$ for $||i-j||\geqslant 2,$\\
$\sigma_0, \ \sigma_1,\ \sigma_2,  \dots , \sigma_{n-1}$ &$\sigma_i \sigma_{i+1}\sigma_i=\sigma_{i+1}\sigma_i \sigma_{i+1}$\\
&$\sigma_{0}=\tau\sigma_{n-1}\tau^{-1}$
\end{tabular}\right\rangle, $\\

where $\tau$ is defined in Lemma 2.1, $i, \ j \in \mathbb{Z}_n,$ and the norm of  $m$ in $\mathbb{Z}_n$ is defined by $||m||=\min\{m, n-m\}$ (here we identify the element $m$ of $\mathbb{Z}_n$ with a number $m\in\{0, \ 1, \dots, n-1\}\subseteq \mathbb{Z}$). \\

 Let $\rho: B_n\to GL_r(\mathbb{C})$ a matrix representation of $B_n$ of \\ dimension $r.$\\
 
 Throughout the paper, we are going to use the following notations.
 
 Let $\rho(\sigma_i)=C_i=I+A_i, \ i\in {\mathbb{Z}_n},$
 
$\rho(\tau)=T$

$\rho(\Delta)=D.$\\
  
 By using braid relations and Lemma 2.1, we obtain the following: \\

{\bf Lemma 2.2.} Let $\rho: B_n\to GL_r(\mathbb{C})$ a matrix representation of $B_n$. Then:\\ 

(a) $TC_{i-1}T^{-1}=C_{i},$ $ \  i\in \mathbb{Z}_n;$ 

(b) $C_iT^{m}=T^{m}C_{i-m} ,\  i, m \in \mathbb{Z}_n;$
      
(c) $TA_iT^{-1}=A_{i+1},$ $i\in {\mathbb{Z}_n};$

(d) $DA_iD^{-1}=A_{n-i},$ $i=1, \ 2, \dots , n-1;$

 (e) For $||i-j||\geqslant 2,$   $A_iA_j=A_jA_i,$ $ i,\ j \ \in {\mathbb{Z}_n};$

(f) $A_i+A_i^2+A_iA_{i+1}A_i=A_{i+1}+A_{i+1}^2+A_{i+1}A_{i}A_{i+1}, \ \  i  \in {\mathbb{Z}_n}.$

(g)  If $\rho: B_n\to GL_r(\mathbb{C})$ is {\it irreducible}    representation of $B_n,$ then\\ $\rho(\tau^n)=\rho(\Delta^2)=T^n=D^2=\alpha I$ for some $\alpha\in \mathbb{C}^{*}$\\

{\it Proof of }(g). It follows from the fact that $\tau^n=\Delta^2$ ( which represents the positive  full twist) is a central element in $B_n,$ and $\rho$ is irreducible. \begin{flushright}$\square$\end{flushright}
 
{\bf Definition 2.3.} Let $\rho: B_n\to GL_r(\mathbb{C})$ be a matrix representation of $B_n$ with $\rho(\sigma_i)=I+A_i, \ i\in {\mathbb{Z}_n}.$\\
 Define the {\bf corank} of a representation $\rho$  by $corank(\rho)=rk(A_1).$\\

By Lemma 2.2(c), we have $corank(\rho)=rk(A_i)\ \ \forall i\in {\mathbb{Z}_n}.$\\

If $\rho$ is a one-dimensional representation of $B_n,$ it is constant, and we will denote it by 
$\chi (y):B_n \to {\mathbb{C^{*}}},$ where $\chi (y)(\sigma_i)=y,$ $ y\in {\mathbb{C}^{*}}$\\ for all $  i \in \mathbb{Z}_n.$

\section{The Reduction Theorem}

  In this section we will prove {\bf The Reduction Theorem} (Theorem 3.8), which will allow us to investigate the irreducible representation of $B_n$ of dimension $r\leqslant 2n-9$ by dealing only with the irreducible representations of relatively small corank. Since the complete classification of the irreducible representations of $B_n$ of dimension $r$ for $r\leqslant n$ was given in \cite{Formanek}, \cite{S},\cite{Lee} and  \cite{FLSV}, we will formulate and prove all statements in this section for the dimension $r\geqslant n+1,$ even though some of them hold for smaller values of $r$ as well. The restriction $n\geqslant 10$ used in Lemma 3.1,  Theorems 3.3, 3.7, 3.8, and Corollary 3.9, provides $n+1\leqslant 2n-9.$
  
   Throughout this section we will use the notations introduced in section 2.\\

     {\bf  Lemma 3.1.} Let $\rho: B_n\to GL_r(\mathbb{C})$ be an {\it irreducible} matrix representation of $B_n$ of dimension $r ,\ n+1\leqslant r\leqslant 2n-9$ for $n\geqslant 10.$ Then there exists an eigenvalue $\lambda$ of $\rho(\sigma_{n-1})$ such that the {\it  largest } block corresponding to $\lambda$ in the Jordan normal form of $\rho(\sigma_{n-1})$  has multiplicity $d\leqslant n-5.$

      \begin{proof}
    
    Consider the Jordan normal form of $\rho(\sigma_{n-1}).$ Suppose that for every eigenvalue  of $\rho(\sigma_{n-1}),$ the largest block  corresponding to that eigenvalue has multiplicity $d\geqslant n-4.$\\
    
     If $\rho(\sigma_{n-1})$ has two or more distinct eigenvalues, then the dimension \\$r\geqslant (n-4)+(n-4)=2n-8,$ a contradiction with  $r\leqslant 2n-9.$\\
    
     Thus,  $\rho(\sigma_{n-1})$ has exactly one eigenvalue $\lambda.$ For this eigenvalue, the largest block  must be a $1\times 1$ block. Indeed, if the largest block with the multiplicity $d\geqslant n-4$ has size $2\times 2$ or larger, then the dimension $r\geqslant 2(n-4)=2n-8,$ and again we get a contradiction with $r\leqslant 2n-9.$\\
     
     So,  $\rho(\sigma_{n-1})$ has exactly one eigenvalue, and each block of the Jordan normal form of $\rho(\sigma_{n-1})$ is a $1\times 1$ block, thus, $\rho(\sigma_{n-1})=\lambda I.$  Due to the conjugation in $B_n,$ each of $\rho(\sigma_{i})$ has the same Jordan normal form for all $i=1,\ 2,\dots , n-1.$\\
     
     So, $\rho(\sigma_{i})=\lambda I$ for all $ i=1, \ 2, \dots , n-1, $  a contradiction with irreducibility of $\rho.$
  \end{proof}

     To prove the next theorem, we will  recall the lemma from \cite{S}, based on the results of Formanek \cite{Formanek} that $B_{n-2}$ has no irreducible representations of dimensions between 2 and $n-5$ (inclusive) for $n$ large enough (see \cite{S} for details).\\
     
      {\bf Lemma 3.2} ( \cite{S}, Lemma 6.4.) Let $\rho: B_n\to GL_r(\mathbb{C})$ be a matrix representation of $B_n$ of dimension $r ,$ where $n\geqslant 6.$ Suppose that $\lambda $ is an eigenvalue of $\rho(\sigma_{n-1}),$ and the largest Jordan
block corresponding to $\lambda $ has multiplicity $d.$

If $d\leq n-5,$ then there exists a one-dimensional subspace of $\mathbb{C}^{r},$ invariant under $B_{n-2}\times <\sigma_{n-1}>=<\sigma_1, \dots , \sigma_{n-3}>\times <\sigma_{n-1}>.$  \\
     
   {\bf  Theorem 3.3.} Let $\rho: B_n\to GL_r(\mathbb{C})$ be an {\it irreducible} matrix representation of $B_n$ of dimension $r ,\ n+1\leqslant r\leqslant 2n-9,$ $n\geqslant 10.$   Then there exists a one-dimensional subspace of $\mathbb{C}^{r}$ invariant under
   ${B_{n-2}\times <\sigma_{n-1}>}.$\\

 \begin{proof} By Lemma 3.1, there exists an eigenvalue $\lambda$ of $\rho(\sigma_{n-1}),$ such that  the largest  Jordan block of $\rho(\sigma_{n-1})$ corresponding to $\lambda$ has multiplicity $d\leqslant n-5.$ Then by Lemma 3.2,  $\mathbb{C}^{r}$ has a one-dimensional subspace invariant under
   ${B_{n-2}\times <\sigma_{n-1}>}.$
   \end{proof}

   {\bf Lemma 3.4. } Let $\rho: B_n\to GL_r(\mathbb{C})$  be an {\it irreducible} matrix \\ representation of $B_n$ of dimension $r\geqslant n+1,$ $n\geqslant 3.$  Let $v\in \mathbb{C}^{r}, \ v\neq 0,$ be a vector such that $L=span \{v \}$ is invariant under $B_{n-2}\times <\sigma_{n-1}>.$ \\ Then the vectors $v, \ Tv, \ T^2v,\  \dots ,T^{n-3}v$ are linearly independent.\\

     \begin{proof}
     
       Since $\rho|_{B_{n-2}\times <\sigma_{n-1}>}:L\to L$ is a one-dimensional representation of $B_{n-2}\times <\sigma_{n-1}>,$ then $\rho(\sigma_1)v=\rho(\sigma_2)v=\dots =\rho(\sigma_{n-3})v=yv$ and $\rho(\sigma_{n-1})v=xv$ for some $x,y\in \mathbb{C}^*,$ or, in our notations,\\
      $C_iv=yv,\ i=1, \ \dots n-3,$ and $C_{n-1}v=xv.$\\
     Then, by Lemma 2.2(b), we have\\ 
     
     $C_{n-2}(T^{m}v)=T^{m}(C_{n-2-m})v=T^{m}(yv)=y(T^mv)$ \\ for $m=1,\ 2, \dots, n-3.$\\
     
     Suppose that the  vectors $v, \ Tv, \ T^2v,\  \dots ,T^{n-3}v$ are linearly dependent.  Let  $ a_0(v)+a_1(Tv)+\dots +a_{n-3 }(T^{n-3}v)=0$ be a non-trivial linear combination of non-zero vectors, and let $i$ be the smallest index, such that the coefficient $a_i\neq 0.$ By left-multiplying the above linear combination by $T^{-i},$ we get another non-trivial linear combination with  a non-zero coefficient for the vector $v.$ Thus, \\$v\in span \{Tv, \ T^2v,\  \dots ,T^{n-3}v\}.$ \\ 
     
     Let  $v=\sum\limits_{i=1}^{n-3}
b_i(T^{i}v).$ Then\\

     $C_{n-2}v=\sum\limits_{i=1}^{n-3}
b_iC_{n-2}(T^{i}v)=\sum\limits_{i=1}^{n-3}b_iy(T^{i}v)=y\sum\limits_{i=1}^{n-3}b_i(T^{i}v)=yv\in L,$\\

 so $C_iv\in span\{v\}=L$ for every $i=1,\ 2, \dots, n-1.$ Thus, the one-dimensional subspace $L$ is invariant under  $B_n,$ which contradicts the irreducibility of $\rho.$\\
  Therefore, the vectors $v, \ Tv, \ T^2v,\  \dots ,T^{n-3}v$ are linearly independent.

     \end{proof}

      {\bf Corollary 3.5.} Let $\rho: B_n\to GL_r(\mathbb{C})$  be an {\it irreducible} matrix \\ representation of $B_n$ of dimension $r\geqslant n+1,$ $n\geqslant 3.$  Suppose \\$v\in \mathbb{C}^{r}, \ v\neq 0,$ is a vector such that $L=span \{v \}$ is invariant under $B_{n-2}\times <\sigma_{n-1}>.$ \\ Then for any fixed $i\in \mathbb{Z}_n,$ the $n-2$ consecutive vectors \\$T^{i}v, \ T^{i+1}v, \ \dots  ,T^{i+n-3}v$ (all powers of $T$ are taken modulo $n$)  are linearly independent.

      \begin{proof}By applying left multiplication by $T^{-i}$ and Lemma 2.2(g), the statement
     immediately follows from Lemma 3.4.

      \end{proof}
   
    {\bf Remark:} From the proof of Lemma 3.4, we can see that we have $n-3$ linearly independent vectors  $Tv, \ T^2v,\  \dots ,T^{n-3}v,$ such that each of them is in $Ker(C_{n-2}-yI).$ That means that \\$dim(Ker(C_{n-2}-yI))\geqslant n-3,$ or, equivalently, \\$dim(Im(C_{i}-yI))=dim(Im(C_{n-2}-yI))\leqslant r-(n-3)=r-n+3.$ However, as shown in Theorem 3.6, the stronger statement is actually true.\\

        {\bf Theorem 3.6. } Let $\rho: B_n\to GL_r(\mathbb{C})$ be an {\it irreducible} matrix representation of $B_n$ of dimension $r,$ $r\geqslant n+1,$ $n\geqslant 5.$  Suppose that $L=span \{v \}$ is a one-dimensional subspace of $\mathbb{C}^{r},$ invariant under\\
   ${B_{n-2}\times <\sigma_{n-1}>}=<\sigma_1, \dots , \sigma_{n-3}>\times <\sigma_{n-1}>.$  \\
   
   Then $\exists\  y\in \mathbb{C}^{*},$ such that $dim (Im(\rho(\sigma_1)-yI)\leqslant r-n+2.$

    \begin{proof} 
    
    1) As in Lemma 3.4,  consider the values  $x, y \in \mathbb{C}^*,$ such that \\ $C_1v=C_2v=\dots =C_{n-3}v=yv$  and $C_{n-1}v=xv.$ \\ 
    
   We are going to show that for this value of $y$\\
   $dim (Im(\rho(\sigma_1)-yI)\leqslant r-n+2.$\\
   
    Suppose not. Then \\$dim (Im(\rho(\sigma_1)-yI)=dim (Im(C_1-yI))\geqslant  r-n+3$ and \\$dim (Ker(C_1-yI))=r-dim (Im(C_{1}-yI))\leqslant r-(r-n+3)=n-3.$\\
    
    In addition, since all $C_i$ are conjugated by $T,$\\
    $dim (Ker(C_i-yI))=dim (Ker(C_1-yI))\leqslant n-3$ for $i\in \mathbb{Z}_n.$\\
    
  2)  Let $S_i=Ker(C_i-yI), \ i \in \mathbb{Z}_n.$\\
  
  Consider $S_1=Ker(C_1-yI).$ \\ By Lemma 2.2(b),
      $C_{1}(T^{m}v)=T^{m}(C_{1-m})v=T^{m}(yv)=y(T^mv)$ for $m=4,\ 5, \dots, n-1,\ 0,$ since $C_iv=yv$ for $i=1, \ \dots, n-3,$ thus, \\$span\{ T^4v, T^5v, \dots, T^{n-1}v, v\}\subseteq S_1.$\\
      
      By Corollary 3.5, the $n-3$ vectors $T^4v, T^5v, \dots, T^{n-1}v, v$ are linearly independent  by being a subset of linearly independent vectors  \\$T^4v, T^5v, \dots, T^{n-1}v, v, Tv.$ So,\\
      $n-3=dim(span\{ T^4v, T^5v, \dots, T^{n-1}v, v\})\leqslant dim(S_1)\leqslant n-3,$ and thus 
      $ S_{1}=Ker(C_{1}-yI)=span \{ T^4v, T^5v, \dots, T^{n-1}v, v\}$ with \\$dim(S_1)= n-3.$ \\
      
      Now, for any $k\in \mathbb{Z}_n,$ by  multiplying the equation $C_{1}(T^{m}v)=y(T^mv)$  for $m=4,\ 5, \dots, n-1,\ 0,$ by $T^k,$ and using Lemma 2.2(b), we get \\$T^kC_{1}(T^{m}v)=C_{1+k}(T^{m+k}v)=y(T^{k+m}v).$ By using Corollary 3.5 and Lemma 2.2(g), similarly to the above argument, we obtain \\
    $$ S_{i}=Ker(C_{i}-yI)=span \{ T^{i+3}v, T^{i+4}v, \dots T^{n-1}v, v, \dots, T^{i-1}v\}$$     for $i\in \mathbb{Z}_n.$\\

   3)  Consider  $S=span\{ v, Tv, T^2v, \dots, T^{n-1}v\}.$ This is a non-trivial subspace, invariant under $T$  with   $dim \ S \leqslant n.$\\ In addition, $S_i\subseteq S$ for all $i\in \mathbb{Z}_n.$ \\
   
   We claim that $S$ is invariant under $B_n.$\\
    
 First, let's show that $S$ is invariant under $C_1.$ \\
         For $m=4,\ 5, \dots, n-1, \ 0,$ we have  $T^mv\in S_1,$ so \\$C_1(T^mv)=y(T^mv)\in S,$ and it remains to show that  $C_1(Tv), \ C_1(T^2v),$ and $C_1(T^3v) $ are in $S.$ Indeed,\\
     
      $C_1(T^2v)=T^2C_{n-1}v=T^2(xv)=xT^2v\in S.$\\
     
     For $n\geqslant 5,$ $Tv\in S_3\Longrightarrow (C_3-yI)(Tv)=0$ and \\$ 0=C_1[(C_3-yI)(Tv)]=(C_3-yI)(C_1(Tv))\Longrightarrow C_1(Tv)\in S_3\subseteq S.$\\
     
     Similarly,\\ $T^3v\in S_{n-1} \Longrightarrow (C_{n-1}-yI)(C_1(T^3v))=0 \Longrightarrow C_1(T^3v)\in S_{n-1}\subseteq S.$\\
     
    Now, let's show that $S$ is invariant under $C_j$ for every $j=2, \dots n-1.$\\ For $m\in \mathbb{Z}_n,$ we have\\  
      $ C_j(T^mv)=T^{j-1}C_{j-(j-1)}(T^{m-(j-1)}v)=T^{j-1}C_{1}(T^{m-j+1}v)\in S,$\\ since $S$ is invariant under  $C_1,$ and $S$ is invariant under $T.$\\

     Thus, we have a proper subspace $S$ invariant under $B_n,$ a contradiction with  the irreducibility of $\rho.$ So, for this value of $y$\\$dim (Im(\rho(\sigma_1)-yI)\leq r-n+2.$

    \end{proof}
     
       {\bf Theorem 3.7.} Let $\rho: B_n\to GL_r(\mathbb{C})$ be an {\it irreducible} matrix \\ representation of $B_n$ of dimension $r ,\ n+1\leqslant r\leqslant 2n-9,$ $n\geqslant 10.$ Then $\exists\  y\in \mathbb{C}^{*},$ such that $dim (Im(\rho(\sigma_1)-yI)\leqslant r-n+2.$
   
  \begin{proof} 
  
  By Theorem 3.3,  there exists a one-dimensional subspace of $\mathbb{C}^{r}$ invariant under
   ${B_{n-2}\times <\sigma_{n-1}>}.$ Then, by Theorem 3.6, $\exists\  y\in \mathbb{C}^{*},$ such that $dim (Im(\rho(\sigma_1)-yI)\leqslant r-n+2.$

    \end{proof}

  We will reformulate the above theorem in terms of coranks.
   
     \vskip 0.3cm 
   {\bf Theorem 3.8 (The Reduction Theorem). }\\ Let $\rho: B_n\to GL_r(\mathbb{C})$ be an {\it irreducible} matrix representation of $B_n$ of dimension $r ,\ n+1\leqslant r\leqslant 2n-9,$ $n\geqslant 10.$ Then $\rho$ is equivalent to a tensor product of a one-dimensional representation  and an irreducible representation $\widehat{\rho}$ of dimension $r$ and corank $k,$ where
   $3\leqslant k \leqslant r-n+2.$
   
   \begin{proof}
   
    By  Theorem 3.7, there exists $ y\in \mathbb{C}^{*},$ such that \\$dim (Im(\rho(\sigma_1)-yI)\leqslant r-n+2.$ \\Consider  one-dimensional representations $\chi(y)$ and $\chi(y^{-1}).$ \\Then for every group generator $\sigma_i\in B_n$\\ $\rho(\sigma_i)-yI=\left(\chi(y)\otimes[\chi(y^{-1})\otimes \rho-I]\right)(\sigma_i)$
   
   Since $\rho$ is irreducible, the representation $\widehat{\rho}=\chi(y^{-1})\otimes \rho$ is also irreducible of the same dimension $r,$ and since $y\neq 0,$ \\
  $ corank (\widehat{ \rho})=dim (Im[(\chi(y^{-1})\otimes \rho-I)(\sigma_1)])=dim (Im(\rho(\sigma_1)-yI))\leqslant$\\$\leqslant r-n+2.$
   
    By \cite{Formanek},  Theorem 10, every irreducible representation of corank \\$k=1$ has dimension $r\leqslant n-1.$ By \cite{S}, Theorem 5.5, every irreducible representation of corank $k=2$ has dimension $r=n.$ 
   Since $r\geqslant n+1,$ we have that $3\leqslant corank (\widehat{\rho})\leqslant r-n+2.$

    \end{proof}

 {\bf Corollary 3.9.} Let $\rho: B_n\to GL_r(\mathbb{C})$ be an {\it irreducible} matrix \\ representation of $B_n$ of dimension $r ,\ n+1\leqslant r\leqslant 2n-9,$ $n\geqslant 10.$ Then $\rho$ is equivalent to a tensor product of a one-dimensional representation  and an irreducible representation of dimension $r$ and corank $k,$ where
   $3\leqslant k \leqslant n-7.$

  \begin{proof}  $r\leqslant 2n-9\Longrightarrow k \leqslant 2n-9-n+2=n-7.$
      
     \end{proof}

\section{Friendship Graphs}

 In this section we will describe the graphs associated with representations of $B_n,$ which will help us to investigate the representations of small corank. These graphs were first introduced in \cite{S}. We will quickly review the definitions and some of the results published in \cite{S}, as well as prove some new results which will be used in section 5 to classify irreducible representations of $B_n$ of  corank 3.\\
 
 Let $\rho:B_n\to GL_{r}(\mathbb{C})$ be a matrix representation of $B_n$ where \\ $\rho(\sigma_i)=I+A_i, \ i\in \mathbb{Z}_n.$\\

Each  graph associated with a representation is a finite {\it simple-edged} graph (there is at most one unoriented edge joining two vertices, and  no edge joins a vertex to itself) such that each vertex of a graph  corresponds to an image of a braid group generator.
We will (slightly abusing the notations) denote each vertex correspoding to $\rho(\sigma_i)$  by $A_i.$\\

{\bf Definition 4.1.} Let $A_i$ and $A_j$ be two {\it distinct} vertices, $i, \ j\in\mathbb{Z}_n.$\\

(i)  $A_i$ and $A_{j}, \ i, j \in {\mathbb{Z}_n},$  are {\bf neighbors} if $||i-j||=1.$\\

(ii)  $A_i$ and $A_j$ are {\bf friends,} if $Im(A_i)\cap Im(A_j)\neq \{0\}.$\\

(iii)  $A_i$ and $A_j$ are {\bf true friends,} if \\

$\left\{\begin{array}{l}
A_i+A_i^2+A_iA_{j}A_i=A_{j}+A_{j}^2+A_{j}A_{i}A_{j}\neq 0, {\textrm { if }} ||i-j||=1;\\
\\
A_iA_j=A_jA_i\neq 0 , {\textrm { if }}   ||i-j||\geqslant 2.
\end{array}  \right.$\\
\vskip .5cm

{\bf Definition 4.2.} {\bf The full friendship graph} associated with the representation
$\rho: B_n\to GL_r(\mathbb{C})$  with $\rho(\sigma_i)=I+A_i, \ i\in {\mathbb{Z}_n},$ is the simple-edged graph with $ n$ vertices
$A_0 ,\  A_1,\ \dots, A_{n-1}$ such that $A_i$ and $A_j$ are connected by an edge  if and only if $A_i$ and $A_j$ are friends.\\

{\bf Definition 4.3.} For two distinct vertices $A$ and $B,$  define:\\

(i)  $f(A,B)=dim\left(Im(A)\cap Im(B)\right);$ \\

(ii) $tf(A,B)=\left\{\begin{array}{l}
dim(Im\left(A+A^2+ABA\right)), {\textrm { if $A$ and $B$ are neighbors}}\\
\\
dim(Im\left(AB\right)), {\textrm { if $A$ and $B$ are not neighbors}}
\end{array}  \right.$
\vskip .5cm

{\bf Lemma 4.4.}  (i) $f(A,B)=f(B,A);$ \\

(ii) $tf(A,B)=tf(B,A).$

\begin{proof}

 Both statements easily follow from the definitions of $f(A,B)$ and $tf(A,B).$
 
 \end{proof}
    
In the above notations, we can say that the vertices $A$ and $B$ are connected by an edge in a  full friendship graph if and only if $f(A,B)>0.$\\

One can also consider the {\bf friendship graph} obtained by removing the extra vertex $A_0$ and all edges incident to it from the full friendship graph (see \cite{S} for details). In addition, in some cases it might be useful to consider {\bf true friendship graphs}, defined in a similar manner to the frienship graphs,  as well as {\bf weighted} friendship and true friendship graphs, assigning the edges' weights to be $f(A,B)$ or $tf(A,B)$ respectively. We will not go into details of these explorations in this paper. Our main interest here will be the full friendship graph, and we will refer to it simply as ``friendship graph'' for the remainder of the paper. \\

The following lemma is the strengthening of Lemma 3.1 from \cite{S}, reformulated in the above notations.\\
 
{\bf Lemma 4.5.} For any two distinct vertices $A$ and $B,$  \\$f(A,B)\geq tf(A,B).$

\begin{proof}

 1) For $A$ and $B$  not neighbors, we have $AB = BA$ so,\\
$Im(A) \cap  Im(B) \supseteq Im(AB) \cap Im(BA) = Im(AB) $, and \\
$f(A,B)\geq tf(A,B).$\\

2) For  $A$ and $B$  neighbors, \\
$A(I + A + BA) = A + A^2+ ABA = B + B^2+ BAB = B(I + B + AB),$
and \\ 
$Im(A)\cap Im(B) \supseteq Im(A(I + A + BA))\cap (Im(B(I + B + AB))=$

$=Im (A + A^2+ ABA ) ,$ so 
$f(A,B)\geq tf(A,B).$

\end{proof}
    
    The next lemma is instrumental in the classification of the possible  friendship graphs for a  representation of $B_n.$\\
    
      {\bf Lemma 4.6.} {\bf (Enhanced Lemma about Friends)} For any three distinct vertices $A,$ $B,$ and $C, $ such that  $A$ and $B$ are neighbors, and $A$ and $C$ are  not neighbors, 
      $$f(B,C)\leq tf(A,B)+tf(A,C)$$
  
  \begin{proof}
  
   Consider a linear map 
   
   $\varphi: Im(B) \cap Im(C)\mapsto \mathbb{C}^r \bigoplus \mathbb{C}^r$ given by 
 
   $\varphi (z) = ((I+B+BA)z,Az)$
     \vskip 0.3cm
    1)  $\forall z\in  Im(B) \cap Im(C)\ \exists x\ ,y\in \mathbb{C}^r$ such that $z=Bx=Cy$ and
   
     $\varphi (z) = ((I+B+BA)z,Az)=$
     
     $=((B+B^2+BAB)x,ACy)\in Im(B+B^2+BAB) \bigoplus Im(AC).$
     
     Thus, $Im(\varphi)\subseteq Im(B+B^2+BAB) \bigoplus Im(AC).$
       \vskip 0.3cm
     2) $Ker(\varphi)=\{0\}$  (so $\varphi$ is injective). Indeed, if
     
     $\varphi(z) = 0$ then $(I+B+BA)z=Az = 0 \Longrightarrow (I+B)z=0$ and since $I+B$ is invertible, $z=0.$
 
       \vskip 0.3cm
     3) From 1) and 2) it follows that $$dim(Im(B) \cap Im(C))\leq dim (Im(B+B^2+BAB))+dim(Im(AC)),$$
     or, in terms of $f$ and $tf,$ $$f(B,C)\leq tf(A,B)+tf(A,C)$$
     
     \end{proof}
      
   {\bf Remark.} The Lemma About Friends (Lemma  3.3 from \cite{S}) is an easy consequence of the Lemma 4.6. \\Indeed,  if $f(A,B)=0$ and $f(B,C)\neq 0,$ then \\
   $0\leqslant tf(A,B)\leqslant f(A,B)=0,$ and \\
    $0\neq f(B,C)\leq  tf(A,B)+tf(A,C)=0+tf(A,C),$ so \\$tf(A,C)\neq 0.$\\

 The following lemma will allow us to talk about action of $\mathbb{Z}_n$ on the friendship graph by cyclically permuting the vertices.\\
 
{\bf Lemma 4.7.} For any $i,\  j, \ k  \in {\mathbb{Z}_n},\ i\neq j,$ \\

(i) $f(A_i,A_j)=f(A_{i+k},A_{j+k});$\\

(ii) $tf(A_i,A_j)=tf(A_{i+k},A_{j+k}).$

\begin{proof}

  Both statements easily follow from the fact that $A_{i+k}=T^kA_iT^{-k}$ for all $i, \ k \ \in \mathbb{Z}_n.$
  
\end{proof}

{\bf Definition 4.8.}  For all $k=1, \ 2, \ \dots , n-1,$ define \\$f(k)=f(A_1,A_{1+k})$ and $tf(k)=tf(A_1,A_{1+k}).$\\

{\bf Remark.}  From the definition and  Lemma 4.7,  one can see that both $f(k)$ and $tf(k)$ can be viewed as the functions describing the corresponding dimensions for the vertices that are $k$ vertices apart from each other, that is $f(A_i,A_j)=f(|i-j|)$ and $tf(A_i,A_j)=tf(|i-j|)$ for all $  i\neq j$ (here $|i-j|$ is the distance between the natural numbers $i$ and $j$).\\

 Clearly, $f(k)=f(n-k)$ and $tf(k)=tf(n-k)$ for all \\$k=1, \ 2, \dots , n-1,$ or, equivalently,\\ $f(A_i,A_j)=f(||i-j||)$ and $tf(A_i,A_j)=tf(||i-j||)$ for all $  i\neq j.$\\

Now, let's  recast two theorems from \cite{S}.\\

  {\bf Theorem 4.9.}(\cite{S}, Theorem 3.4.) Let $\rho: B_n\to GL_r(\mathbb{C})$ be a matrix representation of $B_n$ of dimension $r.$ Then one of the
following holds:\\
 (a) The full friendship graph is totally disconnected (no friends at all);\\
(b) The full friendship graph has an edge between  $A_i$ and $A_{i+1}$ \\for all $i;$\\
(c) The full friendship graph has an edge between  $A_i$ and  $A_j$ whenever
$A_i$ and  $A_j$ are not neighbors.\\

 {\bf Theorem 4.10.} (\cite{S}, Theorem 3.8.) Let $\rho: B_n\to GL_r(\mathbb{C})$ be an {\it irreducible} matrix representation of $B_n$ of dimension $r$ with  a totally disconnected associated friendship graph.\\
  Then $r\leqslant n-1.$\\
  
  We will formulate the corollary of the Theorems 4.9 and 4.10 in terms of  $f$ and $tf.$\\
  
  {\bf Corollary 4.11.} Let $\rho: B_n\to GL_r(\mathbb{C})$ be an {\it irreducible} matrix representation of $B_n$ of dimension $r\geqslant n+1,$ $n\geqslant 3.$ Then one of the
following holds:\\

(a) $f(1)\geqslant 1.$

(b) $f(1)=f(n-1)=0$ and $f(k)\geqslant 1\ \forall k=2,\dots n-2.$\\

Now we will establish some properties that we will use in the next section.\\ 

  {\bf Lemma 4.12.}  Let $\rho: B_n\to GL_r(\mathbb{C})$ be a matrix representation \\of $B_n.$ \\
      
      Then  $tf(2)=tf(3)=\dots =tf(n-2).$
      
      \begin{proof}
      
       To prove this statement, we are going to show that\\
        $tf(A_1, \ A_k)=tf(A_1,\  A_j)$ for $k, \ j = 3,\  4, \dots , n-1.$\\
       
       For any $k=3, \ \dots , n-1,$ \\
       $A_kA_1=A_1A_k,$ so $\forall w \in Im(A_1) \Rightarrow A_kw\in Im(A_1),$ so $A_k$ acts \\ on $Im(A_1).$ Denote by $\widetilde{A}_k=A_k|_{Im(A_1)}.$ Then\\
       
       $tf(A_1, \ A_k)=dim(Im(A_kA_1))=dim (Im(\widetilde{A}_k)),$ $k=3, 4, \dots, n-1.$ \\
       
       Since $A_k$ and $A_j$ are conjugated in $B_{n-2}=<\sigma_3, \ \sigma_4, \ \dots , \sigma_{n-1}>,$ \\we have  $dim (Im(\widetilde{A}_k))=dim (Im(\widetilde{A}_j)), $  so \\
    $tf(A_1, \ A_k)=tf(A_1,\  A_j)$ for $k, \ j = 3,\  4, \dots , n-1.$\\
  
  \end{proof}

      {\bf Lemma 4.13.}  Let $\rho: B_n\to GL_r(\mathbb{C})$ be a matrix representation \\ of $B_n$ with 
      $f(1)=0.$ Then\\
      $f(2)=tf(2)=f(3)=tf(3)=\dots =f(n-2)=tf(n-2).$

     \begin{proof}
     
      By  definition and Lemma 4.7, $f(2)=f(n-2).$ \\
      
      By Lemma 4.5, $tf(k)\leqslant f(k)$ for all $k=1, \ 2, \ \dots , n-1.$\\
      
     By applying Enhanced Lemma about Friends (Lemma 4.6) to \\$B=A_{k}, \ A=A_{k+1}$ and $C=A_1$ for $k=3, \dots , n-2,$ and using \\$0\leqslant tf(A_k,A_{k+1})\leqslant f(A_k,A_{k+1})=0,$ we get\\
     
     $f(A_{1},A_{k})\leqslant tf(A_{k},A_{k+1})+tf(A_{1},A_{k+1})\Longrightarrow$\\
     
     $f(k-1)\leqslant 0+tf(k)=tf(k)\leqslant f(k)$ for all $k=3, \dots n-2.$\\

      Thus,  $f(2)\leqslant tf(3)\leqslant f(3)\leqslant \dots \leqslant tf(n-2)\leqslant f(n-2)=f(2),$ \\which gives the required statement.
      
  \end{proof}
        
  {\bf Lemma 4.14.} Let $\rho: B_n\to GL_r(\mathbb{C})$ be an {\it irreducible} matrix representation of $B_n$ of dimension $r\geqslant n+1,$ $n\geqslant 5$ with \\ $corank(\rho)=dim(Im(A_1))=k$  where $3\leqslant k\leqslant r-1.$ \\ Then $f(j)<k$ for $j=1, \ 2, \dots , n-1.$
 
 \begin{proof}
   1) Suppose that \\
  $f(j-1)=dim\left(Im(A_1)\cap (Im(A_{j}) \right)=k=dim(Im(A_1))$ \\ for some $j=2, \ 3, \dots , n.$ Then, since $ \left(Im(A_1)\cap (Im(A_{j}) \right)\subseteq Im(A_1),$ it follows that $Im(A_1)\cap (Im(A_{j})=Im(A_1).$\\
  
  Similarly, $ \left(Im(A_1)\cap (Im(A_{j}) \right)\subseteq Im(A_j),$ and \\$Im(A_1)\cap (Im(A_{j})=Im(A_j).$ Thus, $Im(A_1)=Im(A_j).$\\
  
  2) Let $W=Im(A_1)=Im(A_j).$  We claim that $W$ is invariant \\under $B_n.$ Indeed, $W$ is invariant under $C_1$ and $C_j.$ We have 2 cases:\\
  
  (a)  $j\neq 3.$ For every $m=3, \dots , n-1, \  m\neq j,$ we have $|m-1|\geqslant 2$ and hence, $C_mw=C_mA_1u=A_1C_mu\in W\  \forall w=A_1u\in Im (A_1)=W.$ \\ And for $m=2,$ we have $|m-j|\geqslant 2,$ and \\$C_mw=C_mA_ju=A_jC_mu\in W\  \forall w=A_ju\in Im (A_j)=W.$ \\
  
  (b) $j=3.$  In this case $W$ is invariant under all $C_m$ for $m\geqslant 4,$ and we only need to check that $W$ is invariant under $C_2.$ \\
  
  Since $n\geqslant 5,$ we  have that \\$f(2)=dim (Im(A_3)\cap Im(A_5))=k,$ (if $n=5$, then use $A_0$ instead of $A_5$) and, hence, $W=Im(A_1)=Im(A_3)=Im(A_5),$ and then \\ $C_2w=C_2A_5u=A_5C_2u \in W\  \forall w=A_5u\in Im (A_5)=W.$\\
  
  3) For $3\leqslant k\leqslant r-1,$ we have that $W$ is a proper invariant subspace of $\mathbb{C}^r,$ which contradicts with the irreducibility of $\rho.$
  
  \end{proof}
 
  In those cases, when the edges of the friendship graph represent  one-dimensional subspaces, it is convenient to consider vectors that span these subspaces. We are going to set the notations that we will use in section 5 to investigate the friendship graphs corresponding to our representations.\\
    
    {\bf Definition 4.15.}  Let $\rho: B_n\to GL_r(\mathbb{C})$ be a matrix representation of $B_n,$ $n\geqslant 3.$ Suppose that for some $i\neq j, \ i, \ j\in \mathbb{Z}_n,$ \\$dim(Im(A_i))\cap Im(A_j))=1,$ and suppose $v\neq 0$ is a vector such that $span \{v\}=Im(A_i)\cap Im(A_{j}).$\\
    
    1) Set $v_{i,j}=v_{j,i}=v.$  \\

2)  If $||i-j||\geqslant 2,$ we will call the vector $v_{i,j}$ {\it a diagonal of the full friendship graph } or, for short, {\it a diagonal.}\\
  In addition we will say that the diagonal $v_{i,j}$ is {\it coming out of the vertex $A_i$,} or, for short, {\it coming out of $A_i$} (as well as $v_{i,j}$ is coming out of $A_j$).\\
  We also will say that the diagonal $v_{i,j}$ {\it connects $A_i$ and $A_j.$}\\

{\bf Lemma 4.16.} Let $\rho: B_n\to GL_r(\mathbb{C})$ be a matrix representation \\of $B_n,$ $n\geqslant 3.$ Suppose that 3 distinct vertices $A_i, \ A_j$ and $A_k, $ \\$ i, \ j, \ k\in \mathbb{Z}_n$ are pairwise connected by  edges representing one-dimensional subspaces. Namely, $dim(Im(A_i))\cap Im(A_j))=dim(Im(A_i))\cap Im(A_k))=$\\$=dim(Im(A_j))\cap Im(A_k))=1.$\\

If $Im(A_i)\cap Im(A_j)=Im(A_i)\cap Im(A_k)$ then \\$Im(A_j)\cap Im(A_k)=Im(A_i)\cap Im(A_j)=Im(A_i)\cap Im(A_k).$ \\
\vskip -.3cm
\setlength{\unitlength}{0.5cm}
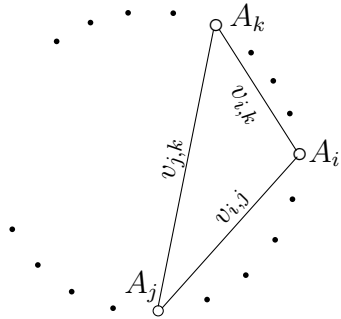
\begin{picture}(10,10)(-5,-5)
\put(3.85,0){$\circ$}
\put(4.3,0){$A_{i}$}

\put(-0.6,-3.5){$A_{j}$}
\put(1.63,3.42){$\circ$}
\put(2.2,3.6){$A_{k}$}
\put(0.3,-3.97){\circle{0.27}}

\put(3.36,2.15){\circle*{0.15}}
\put(3.8,1.28){\circle*{0.15}}
\put(2.77,2.9){\circle*{0.15}}

\put(3.975,0.295){\line(-1,1.57){2.06}}
\put(0.4,-3.9){\line(1,1.12){3.6}}
\put(0.3,-3.85){\line(1,5){1.47}}

\put(0.7,3.95){\circle*{0.15}}
\put(-0.5,3.96){\circle*{0.15}}
\put(-1.5,3.7){\circle*{0.15}}
\put(-2.4,3.2){\circle*{0.15}}

\put(-0.9,-3.9){\circle*{0.15}}
\put(-2,-3.46){\circle*{0.15}}
\put(-2.9,-2.75){\circle*{0.15}}
\put(-3.58,-1.8){\circle*{0.15}}

\put(2,1.9){{\rotatebox{303}{\small{$v_{i,k}$}}}}
\put(0.4,-0.5){{\rotatebox{85}{\small{$v_{j,k}$}}}}
\put(1.8,-1.7){{\rotatebox{60}{\small{$v_{i,j}$}}}}

\put(3.87,-1){\circle*{0.15}}
\put(3.41,-2.08){\circle*{0.15}}
\put(2.63,-3.02){\circle*{0.15}}
\put(1.6,-3.67){\circle*{0.15}}

%\put(0,0){\circle{8}}
%\put(0,0){\circle*{0.08}}

\end{picture}

\vskip -0.3cm

\begin{proof}

 Clearly, if $v_{i,j}\in span\{v_{i,k}\}$ then $v_{i,j}\in  Im(A_j)\cap Im(A_k),$ and the statement follows from the fact that all the subspaces in question are one-dimensional.

\end{proof}
  
 {\bf Cyclicity and Symmetry Arguments}. By Lemma 2.2, parts (c) and (d), we have:\\
 
 (1) $TA_iT^{-1}=A_{i+1},$ $i\in {\mathbb{Z}_n}$

(2) $DA_iD^{-1}=A_{n-i},$ $i=1, \ 2, \dots , n-1.$\\

By using the fact that for two linear operators $A$ and $B$ on a vector space $V,$ conjugated by an invertible $P$ ($B=PAP^{-1}$), if $v\in Im(A)$ then $Pv\in Im(B)$, and other basic facts from linear algebra, we will refer to (1) as {\bf cyclic argument} and to (2) as {\bf symmetric argument}.\\

 Note that $\Delta \sigma_0\Delta^{-1} \neq \sigma_0,$ so the symmetry argument can not be extended to $i=0.$
 
 \section{Irreducible representations of corank 3.}

  The main goal of this paper is  the classification of the irreducible representations of dimension $n+1.$ 
  First, we are going to reduce this problem  to the classification of the irreducible representations of  small corank.\\

   {\bf Theorem 5.1 }  Let $\rho: B_n\to GL_{n+1}(\mathbb{C})$ be an {\it irreducible} matrix representation of $B_n$ of dimension $n+1,$ $n\geqslant 10.$ Then $\rho$ is equivalent to a tensor product of a one-dimensional representation  and an irreducible representation of dimension $n+1$ and corank $3.$

 \begin{proof}
 
  By the Reduction Theorem (Theorem 3.8),  $\rho$ is equivalent to a tensor product of a one-dimensional representation  and an irreducible representation of dimension $r=n+1$ and corank $k$ with \\$3\leqslant k \leq r-n+2=n+1-n+2=3.$
  
  \end{proof}
  
In this section we are going to prove that for $n\geqslant 10$ the braid group  $B_n$ has no irreducible representations  of dimension $r\geqslant n+1$ of corank 3. \\

  Suppose $\rho: B_n\to GL_{r}(\mathbb{C})$ is an {\it irreducible}  representation of $B_n$ of dimension $r\geqslant n+1,$  $n\geqslant 3.$   By Corollary 4.11, there are only two possibilities:\\

 {\bf Case I.} $Im(A_i)\cap Im(A_{i+1})\neq \{0\}$ for all $i\in \mathbb{Z}_n.$\\

 {\bf Case II.}  $Im(A_i)\cap Im(A_{i+1})= \{0\}$ for all $i\in \mathbb{Z}_n$ and\\ $Im(A_i)\cap Im(A_{j})\neq \{0\}$ for $||i-j||\geqslant 2, \ i, \ j\in \mathbb{Z}_n. $\\

   We will consider these two cases separately in subsections 5.1 and 5.2. We will use the notations introduced in sections 2 and 4.\\

  { \centerline{5.1. \bf{ Case I: $Im(A_i)\cap Im(A_{i+1})\neq \{0\}.$}}}
 
  \vskip 0.3cm 
  
 In this subsection we will prove that there are no irreducible representations of $B_n$ of dimension $r\geqslant n+1$ and corank $3$ such that \\ $Im(A_i)\cap Im(A_{i+1})\neq \{0\}$ for $n\geqslant 9$ (Theorem 5.1.24).  First, we will prove that if $Im(A_i)\cap Im(A_{i+1})\neq \{0\}$ then $Im(A_i)\cap Im(A_{i+1})$ must be one-dimensional for all $i$ (Theorem 5.1.5), and then we will prove that in this case the graph is complete  with $dim(Im(A_i)\cap Im(A_{j}))=1$ for all $i\neq j$ (Theorem 5.1.15). We will finish this subsection by investigating the corresponding friendship graph to prove Theorem  5.1.24.\\
 
 We are going to start with the following  lemma.\\

   {\bf Lemma 5.1.1.}  Let $\rho: B_n\to GL_{r}(\mathbb{C})$ be an {\it irreducible}  matrix representation of $B_n$ of dimension $r\geqslant n+1$  such that $dim(Im(A_i))=k,$ where $ 3\leqslant k\leqslant r-1$ and $n \geqslant 5.$ \\
Suppose that for some $j =2, \dots, n-3$ \\$Im(A_1)\cap Im(A_{2})\cap \dots \cap  Im(A_j) \neq \{0\}.$ \\
 Then for every $i\in \mathbb{Z}_n$  \\
 $Im(A_i)\cap Im(A_{i+1})\cap \dots \cap Im(A_{i+j-1})\neq $\\
 $\neq Im(A_{i+1})\cap Im(A_{i+2})\cap\dots \cap Im(A_{i+j})$\\ (the indices are taken modulo $n$).

  \begin{proof}
  
   Suppose that for some $i\in \mathbb{Z}_n$  \\ $Im(A_i)\cap Im(A_{i+1})\cap \dots \cap Im(A_{i+j-1})=$\\$= Im(A_{i+1})\cap Im(A_{i+2})\cap\dots \cap Im(A_{i+j}).$\\ 
  Then, by the cyclic argument, \\
    $Im(A_1)\cap Im(A_{2})\cap \dots  \cap  Im(A_j)  =$\\$=Im(A_2)\cap Im(A_{3})\cap \dots \cap  Im(A_{j+1})=$\\ $=\dots=$\\ $=Im(A_{0})\cap Im(A_{1})\cap \dots \cap Im(A_{j-1})\neq \{0\}.$  \\
   
  Consider the  subspace $S$ defined by\\
   $S=Im(A_1)\cap Im(A_{2})\cap Im(A_{3})\cap \dots \cap Im(A_{n-1})\cap Im(A_{0}).$\\
   It follows from our assumption that  \\$S=Im(A_{i})\cap Im(A_{i+1})\cap  \dots \cap Im(A_{i+j-1})$ for every $i\in \mathbb{Z}_n.$\\
   
  We will now prove  that the following statements are true:\\

   (1) $S\neq \{0\};$
   
   (2) $dim\ S\leqslant k-1;$
   
   (3) $S$ is $B_n-$invariant.\\
   
 Statement (1)  follows from  the hypothesis of the lemma; \\statement (2)  follows  from the  obvious inclusion $S\subseteq Im(A_1)\cap Im(A_{2})$ and Lemma 4.14.
   
   Let's prove statement (3).  Suppose $x\in S,$ we need to show that   $A_mx\in S$ for all $m\in \mathbb{Z}_n.$   \\
 
   We have: $S=Im(A_{m+2})\cap Im(A_{m+3})\cap\dots \cap Im(A_{m+j+1})$ (the indices are taken modulo $n$). Since $j\leqslant n-3, $ $A_m$ is not a neighbor of each of  $A_{m+2}, \ A_{m+3}, \dots , A_{m+j+1},$ thus, $A_m$ commutes with each of  \\$A_{m+2}, \ A_{m+3}, \dots , A_{m+j+1}.$\\ So, if $x\in Im(A_{m+2})\cap Im(A_{m+3})\cap \dots \cap Im(A_{m+j+1}),$ then \\$A_mx\in Im(A_{m+2})\cap Im(A_{m+3})\cap \dots \cap Im(A_{m+j+1})=S.$\\
   
   Since $r>k,$  it follows from (1), (2), and (3) that $S$ is the proper invariant subspace  of $\mathbb{C}^r,$ which contradicts the irreducibility of $\rho.$
   
\end{proof}

    {\bf Lemma 5.1.2.} Let $\rho: B_n\to GL_{r}(\mathbb{C})$ be a matrix representation \\of $B_n$ with $dim(Im(A_i))=3,$ $n\geqslant 3.$ \\Suppose $dim(Im(A_i)\cap Im(A_{i+1}))=2$ for all $i\in \mathbb{Z}_n.$ \\Then  $Im(A_i)\cap Im(A_{i+1})\cap Im(A_{i+2})\neq\{0\}$ for all  
     $i\in \mathbb{Z}_n.$
       
      \begin{proof}
      
       Using the cyclicity argument, it is enough to show that \\$Im(A_1)\cap Im(A_{2})\cap Im(A_{3})\neq\{0\}.$ \\
     
     Indeed, since $Im(A_1)\cap Im(A_{2})+Im(A_2)\cap Im(A_{3}) \subseteq Im(A_2),$\\
     we obtain \\
     
      $ dim(Im(A_1)\cap Im(A_{2})\cap Im(A_{3}))=$
      
      $=dim(Im(A_1)\cap Im(A_{2}))+dim(Im(A_2)\cap Im(A_{3}))-$
      
      $-dim(Im(A_1)\cap Im(A_{2})+Im(A_2)\cap Im(A_{3}))\geqslant$
      
      $\geqslant 2+2-dim( Im(A_{2}))=1.$

    \end{proof}

     {\bf Lemma 5.1.3.} Let $\rho: B_n\to GL_{r}(\mathbb{C})$ be an {\it irreducible} matrix\\ representation of $B_n$ of dimension $r\geqslant n+1$  with $dim(Im(A_i))=3$ and $n\geqslant 5.$ \\Suppose $dim(Im(A_i)\cap Im(A_{i+1}))=2$ for all $i\in \mathbb{Z}_n.$ \\Then  $dim (Im(A_i)\cap Im(A_{i+1})\cap Im(A_{i+2}))=1$ for all  
     $i\in \mathbb{Z}_n.$

     \begin{proof}
     
      Due to the cyclic argument, is enough to show that \\$dim(Im(A_1)\cap Im(A_{2}) \cap Im(A_{3}))=1.$  \\
       
       By Lemma 5.1.2, $dim(Im(A_1)\cap Im(A_{2})\cap Im(A_{3}))\geqslant 1,$ \\
       and  $dim(Im(A_1)\cap Im(A_{2})\cap Im(A_{3}))\leqslant dim( Im(A_1)\cap Im(A_{2}))\leqslant 2.$ \\

If $dim(Im(A_1)\cap Im(A_{2})\cap Im(A_{3}))=2,$ then \\ $Im(A_1)\cap Im(A_{2})\cap Im(A_{3})=Im(A_1)\cap Im(A_{2})$ and  \\$Im(A_1)\cap Im(A_{2})\cap Im(A_{3})=Im(A_{2})\cap Im(A_{3}).$ Thus, \\$Im(A_1)\cap Im(A_{2})=Im(A_{2})\cap Im(A_{3}),$ a contradiction with \\Lemma 5.1.1 for $j=2$ and $k=3.$ Thus,\\ $dim(Im(A_1)\cap Im(A_{2}) \cap Im(A_{3}))=1.$

\end{proof}

{\bf Lemma 5.1.4.} Let $\rho: B_n\to GL_{r}(\mathbb{C})$ be an {\it irreducible} matrix\\ representation of $B_n$ of dimension $r\geqslant n+1$  with $dim(Im(A_i))=3 $ and $n\geqslant 6.$ \\Suppose $dim(Im(A_i)\cap Im(A_{i+1}))=2$ for all $i\in \mathbb{Z}_n.$ \\For each $i\in \mathbb{Z}_n$  consider non-zero vector $w_i,$ such that \\$span \{w_i\}=Im(A_{i-1})\cap Im(A_{i})\cap Im(A_{i+1}).$\\
Let $W_i=span\{w_{i-1}, \ w_i, \ w_{i+1}\},$ $i\in \mathbb{Z}_n.$\\
Then  $Im (A_i)=W_i$ for all $i\in \mathbb{Z}_n.$

\begin{proof}

 First of all,  by Lemma 5.1.3, $Im(A_{i-1})\cap Im(A_{i})\cap Im(A_{i+1})$ is one-dimensional, so there exists a non-zero vector $w_i$ generating \\ $Im(A_{i-1})\cap Im(A_{i})\cap Im(A_{i+1})$ for each $i\in \mathbb{Z}_n.$\\

Due to the  cyclicity, it is enough to show that the three vectors \\$w_1, \ w_2, \ w_3\in Im(A_2)$ are linearly independent.\\

By Lemma 5.1.1 with $j=3\leqslant n-3$ for $n\geqslant 6,$ we have that \\$Im(A_0)\cap Im(A_{1})\cap Im(A_{2})\neq Im(A_1)\cap Im(A_{2})\cap Im(A_{3}),$ so the \\ vectors $w_1$ and $w_2$ are linearly independent, and $dim(W_2)\geqslant 2.$\\ Similarly, the vectors $w_2$ and $w_3$ are linearly independent.\\

Suppose that  the vectors $w_1, \ w_2, \ w_3$ are linearly dependent. Then \\$dim(W_2)=2$ and  $W_2=span \{ w_1, \ w_2\}=span \{ w_2, \ w_3\}.$\\

Since both $w_1, \ w_2\in  Im(A_1)\cap Im(A_2) $ and \\$dim(span \{ w_1, \ w_2\})=dim (Im(A_1)\cap Im(A_{2}))=2,$ we have $\\W_2=span \{ w_1, \ w_2\}=Im(A_1)\cap Im(A_{2}).$\\
Similarly, $W_2=span \{ w_2, \ w_3\}=Im(A_2)\cap Im(A_{3}),$ so  \\$W_2=Im(A_1)\cap Im(A_{2})=Im(A_2)\cap Im(A_{3}),$ \\a contradiction with Lemma 5.1.1.

\end{proof}

{\bf Theorem 5.1.5.} Let $\rho: B_n\to GL_{r}(\mathbb{C})$ be an {\it irreducible} matrix representation of $B_n$ of dimension $r\geqslant n+1$  with $dim(Im(A_i))=3$ and $Im(A_i)\cap Im(A_{i+1})\neq \{0\}$ for $n\geqslant 6.$ \\Then $dim(Im(A_i)\cap Im(A_{i+1}))=1$ for all $i\in \mathbb{Z}_n.$ 

\begin{proof}

 By Lemma 4.14, $dim(Im(A_i)\cap Im(A_{i+1}))=1$ or \\$dim(Im(A_i)\cap Im(A_{i+1}))=2.$ 

Suppose $dim(Im(A_i)\cap Im(A_{i+1}))=2.$ Then by Lemma 5.1.4, \\  $Im (A_i)=span\{w_{i-1}, \ w_i, \ w_{i+1}\},$ where \\$span \{w_i\}=Im(A_{i-1})\cap Im(A_{i})\cap Im(A_{i+1}),$  $i\in \mathbb{Z}_n.$ Thus, \\$Im(A_1)+Im(A_2)+\dots +Im(A_{n-1})+Im(A_0)=span\{w_1, \ w_2, \dots , w_{n-1}, w_0\}.$ Since $Im(A_1)+Im(A_2)+\dots +Im(A_{n-1})+Im(A_0)$ is a $B_n-$invariant subspace of $\mathbb{C}^r,$ and \\$3\leqslant dim(Im(A_1)+Im(A_2)+\dots +Im(A_{n-1})+Im(A_0))\leqslant n,$ we get a contradiction with the irreducibility of $\rho.$\\

Thus, $dim(Im(A_i)\cap Im(A_{i+1}))=1.$

\end{proof}

  Now, we are going to consider the structure of the friendship graph associated with our representation. First, we are going to prove (Corollary 5.1.8) that in our case we actually have a complete graph, meaning that every two vertices are connected by an edge, and each  edge, in fact, represents a one-dimensional subspace (Theorem 5.1.15).\\
  
  For our convenience, let's introduce the following notations.\\

 For all $i\in\mathbb{Z}_n,$ let $x_i=v_{i,i+1}$ (see Definition 4.15), that is \\$span \{x_i\}=Im(A_i)\cap Im(A_{i+1}),$ and let \\$U=span\{x_0,\ x_1, \ x_2, \ \dots, \ x_{n-1}\}.$\\
  
We will list some of the facts in the following lemma.\\

{\bf Lemma 5.1.6.} Let $\rho: B_n\to GL_{r}(\mathbb{C})$ be an {\it irreducible}  matrix \\ representation of $B_n$ of dimension $r\geqslant n+1$ for $n \geqslant 5$ with  \\
$dim(Im(A_i)\cap Im(A_{i+1}))=1,$ where  $Im(A_i)\cap Im(A_{i+1})=span \{x_i\}.$  \\  
 Then:\\ 
(1) $x_i$ and $x_{i+1}$ are linearly independent for all $i\in \mathbb{Z}_n;$ thus, $dim\ U\geqslant 2;$\\
(2) $Tx_i\in span \{x_{i+1}\}$ for all $i\in \mathbb{Z}_n;$\\
(3) $U\neq \{0\}$ and $dim\ U\leqslant n;$\\
(4) $U$ is invariant under $T.$

\begin{proof}

 (1) follows from Lemma 5.1.1, (2)  follows from the cyclic argument,  (3) is obvious, and (4) follows from (2).
 
\end{proof}

{\bf Lemma 5.1.7.}   Let $\rho: B_n\to GL_{r}(\mathbb{C})$ be an {\it irreducible}  matrix \\ representation of $B_n$ of dimension $r\geqslant  n+1$ with $rk(A_1)=3$ and \\
$dim(Im(A_i)\cap Im(A_{i+1}))=1$  for $n \geqslant 5.$  \\ Then $tf(k)\neq 0$ for every $k= 2, \ 3, \dots , \ n-2.$

\begin{proof}

 Suppose that $tf(k)=0$ for some $k=2, \ 3, \ \dots , n-2.$ Then
by Lemma 4.12, we have $0=tf(k)=tf(2)=tf(3)=\dots=tf(n-2),$ and, in particular, $A_1A_3=A_3A_1=0.$\\

  We claim that in this case $U= span\{x_0, \ x_1, \ \dots , \ x_{n-1}\}$ is $B_n-$invariant.\\

Since $A_kx_i=T^{k-1}A_1T^{-(k-1)}x_i$ and $U$ is $T-$invariant, to show that $U$ is invariant under $B_n,$  it is enough to check that \\$A_1x_i\in U$ for  $i=0,\ 1, \ \dots, {n-1}.$\\

Since $tf(2)=tf(3)=\dots=tf(n-2)=0,$ we have \\$A_1x_2=A_1x_3=\dots =A_1x_{n-1}=0,$ and we only have to show that \\$A_1x_0\in U$ and $A_1x_1\in U.$\\

Let's consider $(I+A_1)x_0.$ Since $tf(2)=0$ and $x_0\in Im(A_0),$ we have $A_2x_0=0.$ Since $x_0\in Im(A_1),$ there exists $z,$ such that $x_0=A_1z.$ Thus,\\

$(I+A_1)x_0=(I+A_1)x_0+0=(I+A_1)x_0+A_1\cdot 0=$

$=(I+A_1)x_0+A_1(A_2x_0)= (I+A_1+A_1A_2)A_1z=$

$=(A_1+A_1^2+A_1A_2A_1)z=(A_2+A_2^2+A_2A_1A_2)z=$

$=A_2(I+A_2+A_1A_2)z\in Im(A_2).$\\

But  $x_0\in Im(A_1),$ so  $(I+A_1)x_0\in Im(A_1),$ so \\
$(I+A_1)x_0\in Im(A_1)\cap Im(A_2)=span \{x_1\}\subseteq U.$ Thus,\\
$A_1x_0=(I+A_1)x_0-x_0\in U.$\\

Similarly, for $x_1=A_1y\in Im(A_1)\cap Im(A_2),$ \\

$(I+A_1)x_1=(I+A_1)x_1+A_1\cdot 0=$

$=(I+A_1)x_1+A_1(A_0x_1)=(I+A_1+A_1A_0)A_1y=$

$=(A_1+A_1^2+A_1A_0A_1)y=(A_0+A_0^2+A_0A_1A_0)y=$

$=A_0(I+A_0+A_1A_0)y\in Im(A_0)\cap Im(A_1)\subseteq U,$\\ 

and hence $A_1x_1\in U.$\\

Thus, $U$ is a $B_n-$invariant proper subspace of $ \mathbb{C}^{r},$ which contradicts the irreducibility of $\rho.$ So, $tf(k)\geqslant 1$ for all $k=2, \ 3, \dots , \ n-2.$ 

\end{proof}

{\bf Corollary 5.1.8.} Let $\rho: B_n\to GL_{r}(\mathbb{C})$ be an {\it irreducible}  matrix representation of $B_n$ of dimension $r\geqslant n+1$ with $rk(A_1)=3$ and \\
$dim(Im(A_i)\cap Im(A_{i+1}))=1$  for $n \geqslant 5.$ \\

Then $Im(A_i)\cap Im(A_j)\neq \{0\}$  for all $i\neq j\in\mathbb{Z}_n.$

\begin{proof} The hypothesis of the lemma states that the statement is true for  $||i-j||=1,$ and \\

 $dim(Im(A_i)\cap Im(A_j))=f(A_i,A_j)\geqslant tf(A_i,A_j)= tf(||i-j||)\geqslant 1$ for all $i,\ j\in \mathbb{Z}_n,$ $||i-j||\geqslant 2 $ by Lemma 4.5 and Lemma 5.1.7.
 
\end{proof}

 Now, let's establish some important properties of the vectors \\$x_0, \ x_1, \dots , x_{n-1}$ and the subspace $U.$\\
     
  {\bf Lemma 5.1.9.}   Let $\rho: B_n\to GL_{r}(\mathbb{C})$ be an {\it irreducible}  matrix \\ representation of $B_n$ of dimension $r\geqslant n+1$  for $n \geqslant 5$ with  \\
$dim(Im(A_i)\cap Im(A_{i+1}))=1,$ where  $Im(A_i)\cap Im(A_{i+1})=span \{x_i\}.$  \\  

Suppose that for some $m,$  the vectors $x_1, \ x_2, \ \dots , x_{m}$ are linearly independent, and   the vectors $x_1, \ x_2, \ \dots , x_{m+1}$ are linearly dependent (indices are taken modulo $n$).\\

Then the vectors $x_1, \ x_2, \ \dots , x_{m}$ form a basis of $U.$  In particular,\\ $dim\ U=m.$  

 \begin{proof}
 
  We will use induction on $k$ to show that \\$x_{m+k}\in span \{x_1, \ x_2, \dots , x_m\}$ for all $k\geqslant 1.$\\
  
   For $k=1,$ the statement is true by the hypothesis of the lemma.\\
   
   Now, $x_{m+k}\in span\{Tx_{m+(k-1)}\},$ and if \\$x_{m+(k-1)}\in span \{x_1, \ x_2, \dots , x_m\},$ then \\$x_{m+k}\in span  \{Tx_1, \ Tx_2, \dots , Tx_m\}= span \{x_2, \ x_3, \dots , x_{m+1}\}$ which, together with $x_{m+1}\in span \{x_1, \ x_2, \dots , x_m\},$ gives \\$x_{m+k}\in  span \{x_1, \ x_2, \dots , x_m\}.$
   
  \end{proof}
  
The following statement is an easy consequence of Lemma 5.1.9 and the cyclic argument.\\ 

{\bf Corollary 5.1.10.} Under conditions of Lemma 5.1.9,  for any fixed $i\in \mathbb{Z}_n,$  $m$ consecutive  vectors $x_i, \ x_{i+1}, \dots , \ x_{i+m-1}$  form a basis of $U$ (the indices are taken modulo $n$).\\

{\bf  Theorem 5.1.11. } Let $\rho: B_n\to GL_r(\mathbb{C})$ be an {\it irreducible} matrix representation of $B_n$ of dimension $r\geqslant n+1 $ for $n\geqslant 5,$ such that  \\
$dim(Im(A_i)\cap Im(A_{i+1}))=1,$ where  $Im(A_i)\cap Im(A_{i+1})=span \{x_i\},$ and 
 $U=span\{x_0, \ x_1,  \ x_2, \dots, \ x_{n-1}\}.$\\

Then $dim\ U\geqslant n-3.$

\begin{proof}

 Since $dim\ U\geqslant 2$, (Lemma 5.1.6, part 1), the statement is trivial for $n=5.$ Suppose now that $n\geqslant 6$ and 
 suppose $m=dim\ U\leqslant n-4.$ We claim that in this case $U$ is invariant under $B_n.$\\ Indeed, for any $i\in \mathbb{Z}_n,$ consider $m$ consecutive vectors starting with $x_{i+2}.$ By Corollary 5.1.10,  they form a basis of $U,$ \\$U=span \{x_{i+2}, \ x_{i+3}, \dots , \ x_{i+m+1}\}$ (the indices are taken modulo $n$). \\Since $m\leqslant n-4, $ then for every  $k=i+2, \dots , i+m+1,$ both \\$A_k$ and $A_{k+1}$ are not neighbors of $A_i,$ and hence, commute with $A_i.$ Thus,  $x_k\in Im(A_k)\cap Im(A_{k+1}),$ so \\$ A_i x_k\in Im(A_k)\cap Im(A_{k+1})=span \{x_k\}\subseteq U.$

Since $U\neq \{0\}$ and  $dim \ U \leqslant n-4<n+1\leqslant r,$ we have that $U$ is a proper  invariant subspace of $\mathbb{C}^r,$  which contradicts the irreducibility of $\rho.$

 \end{proof}

{\bf  Corollary 5.1.12. } Let $\rho: B_n\to GL_r(\mathbb{C})$ be an {\it irreducible} matrix representation of $B_n$ of dimension $r\geqslant n+1 $ for $n\geqslant 5,$ such that  \\
$dim(Im(A_i)\cap Im(A_{i+1}))=1,$ where  $Im(A_i)\cap Im(A_{i+1})=span \{x_i\}.$\\

Then the vectors $x_i, \ x_{i+1}, \dots , \ x_{i+n-4}, \ i\in\mathbb{Z}_n,$ are linearly independent.

\begin{proof}

 Immediately follows from Corollary 5.1.10 and Theorem 5.1.11.
 
 \end{proof}

{\bf  Corollary 5.1.13. } Let $\rho: B_n\to GL_r(\mathbb{C})$ be an {\it irreducible} matrix representation of $B_n$ of dimension $r\geqslant n+1 $ for $n\geqslant 7,$ such that  \\
$dim(Im(A_i)\cap Im(A_{i+1}))=1,$ where  $Im(A_i)\cap Im(A_{i+1})=span \{x_i\}.$\\

Then the vectors $x_i$
  and  $x_j$ are linearly independent for all  $i\neq j,$ \\$ i, \ j \in\mathbb{Z}_n.$

\begin{proof}

 Due to cyclicity, it is enough to show that $x_1$ and $x_j$ are linearly independent for all $j\neq 1, j\in \mathbb{Z}_n.$ By Corollary 5.1.12, for $j\geqslant n-3, $ both vectors $x_1$ and $x_j$ belong to the linearly independent set $x_1, \ x_{2}, \dots , \ x_{n-3},$ and for $j=n-2, \ n-1, \ 0,$ both vectors belong to the linearly independent set $x_5, x_6, \dots, \ x_{0}, x_1$ (since $n\geqslant 7$).
 
\end{proof}

 {\bf  Lemma 5.1.14.} Let $\rho: B_n\to GL_{r}(\mathbb{C})$ be an {\it irreducible}  matrix representation of $B_n$ of dimension $r\geqslant n+1$  with $dim(Im(A_i))=3,$ for $n \geqslant 5,$ such that 
$dim(Im(A_i)\cap Im(A_{i+1}))=1,$ where  \\$Im(A_i)\cap Im(A_{i+1})=span \{x_i\},$ and 
 $U=span\{x_0, \ x_1,  \ x_2, \dots, \ x_{n-1}\}.$
 \\Let $y$ be a non-zero vector $y\in Im(A_1)$ such that $y\notin span \{x_0, \ x_1\}.$\\

Then $y\notin U.$

\begin{proof}

 Since $dim(Im(A_i))=3,$ we have  $Im(A_1)=span \{x_{0}, \ x_1, \ y\}.$ \\By applying the cyclic argument, we have that \\$Im(A_i)=span \{x_{i-1}, \ x_i, \ T^{i-1}y\}$ for all $i\in \mathbb{Z}_n.$\\

Suppose that  $y\in U.$  Then, $Im(A_1) \subseteq U,$ and since $U$ is invariant under $T,$ $Im(A_i) \subseteq U$ for all $i,$ and \\

$dim\ V=dim(Im(A_1)+Im(A_2)+\dots+Im(A_{n-1})+Im(A_0) )\leqslant $

$\leqslant dim\ U \leqslant n, $ a contradiction. 

 \end{proof}
 
 Next, we  prove the theorem about the structure of the friendship graph in the case when the neighbors are connected by an edge.\\
 
   {\bf Theorem 5.1.15.} Let $\rho:B_n\to GL_{r}(\mathbb{C})$ be an {\it irreducible} matrix representation of $B_n$ of dimension $r\geqslant n+1$ with and $rk(A_1)=3,$ where $n\geqslant 9.$\\
Suppose that $Im(A_i)\cap Im(A_{i+1})\neq \{0\}$ for all $i\in \mathbb{Z}_n.$ \\ 

Then $dim(Im(A_i)\cap Im(A_{j}))= 1$ for all $i, j\in \mathbb{Z}_n, \ i\neq j.$

\begin{proof}

 By Theorem 5.1.5, $dim(Im(A_i)\cap Im(A_{j}))= 1$ for $||i-j||=1.$ \\

For $||i-j||\geqslant 2, $ by Corollary 5.1.8, $dim(Im(A_i)\cap Im(A_{j}))\geqslant  1, $ and by Lemma 4.14,  $dim(Im(A_i)\cap Im(A_{j}))\leqslant  2. $ Thus, it remains to show that $dim(Im(A_i)\cap Im(A_{j}))\neq 2.$\\

Due to the cyclic argument, it is enough to prove that \\$dim(Im(A_0)\cap Im(A_{j}))\neq 2$ for $j=2, \ 3, \ \dots , n-2,$ and since \\$f(k)=f(n-k)$ for $k=1, \ \dots, n-1,$   it is enough to prove this statement for $j=2,\  3, \dots, \left[\frac{n}{2}\right],$ where $\left[\frac{n}{2}\right]$ denotes the integer part of $\frac{n}{2}.$\\ Equivalently, it is enough to prove that $dim(Im(A_1)\cap Im(A_{j}))\neq 2$ for all $j=3, \ 4, \dots , \left[\frac{n}{2}\right]+1.$ \\

Suppose that for some $j, \ 3\leqslant j\leqslant \left[\frac{n}{2}\right]+1,$ \\
$dim(Im(A_1)\cap Im(A_{j}))=2.$\\

Consider vectors $x_0, x_1, \ x_2, \dots, \ x_{n-1},$ where  \\$Im(A_i)\cap Im(A_{i+1})=span \{x_i\},$  and a vector $y_1\in Im(A_1)$ such that $y_1\notin span \{x_0, \ x_1\}.$  Then\\ $y_j=T^{j-1}y_1\in Im(A_j),$ \\ $Im(A_1)=span \{y_1, \ x_{0}, \ x_1\},$\\$Im(A_j)=span \{y_j, \ x_{j-1}, \ x_j\}.$\\

$Im(A_1)+Im(A_j)=span \{ x_{0}, \ x_1, \ x_{j-1}, \ x_j, \  y_1, \ y_j\}$ and \\
$dim(Im(A_1)+Im(A_j))=$\\$=dim(Im(A_1))+dim(Im(A_j))-dim(Im(A_1)\cap Im(A_{j}))=$\\$=3+3-2=4.$\\

Since $3\leqslant j\leqslant \left[\frac{n}{2}\right]+1,$ the indices for the vectors  $x_{0}, \ x_1, \ x_{j-1}, \ x_j$ are all distinct. Since $n\geqslant 9$ and  $j\leqslant \left[\frac{n}{2}\right]+1, $ we have \\ 

$n-4-j\geqslant n-4- \left[\frac{n}{2}\right]-1\geqslant n-5-\frac{n}{2}\geqslant -0.5, $\\ and, since $n-4-j\in \mathbb{Z}\Longrightarrow n-4-j\geqslant 0$ or, equivalently, $j\leqslant n-4.$ \\Thus, the set of vectors $x_{0}, \ x_1, \ x_{j-1}, \ x_j$ is a subset of linearly independent vectors $x_{0}, \ x_1,\dots , \ x_{n-5}, \ x_{n-4}$ (Corollary 5.1.12).\\
Thus, $dim(span \{ \ x_{0}, \ x_1, \ x_{j-1}, \ x_j\})=4,$ so, \\$y_1\in span \{ \ x_{0}, \ x_1, \ x_{j-1}, \ x_j\}\subseteq U, $ a contradiction with Lemma 5.1.14.

  \end{proof}

The remainder of this subsection 5.1 is devoted to the representations satisfying the following conditions:\\

\begin{tabular}{cl}
 & Let $\rho: B_n\to GL_{r}(\mathbb{C})$ be an {\it irreducible}  matrix representation \\
 {\bf  ($\bigstar$)}& of $B_n$ of dimension $r\geqslant n+1$ with $rk(A_1)=3$\\
 &and $dim(Im(A_i)\cap Im(A_{j}))=1$ for all $i, \ j \in \mathbb{Z}_n, \ i\neq j.$
\end{tabular}\\

Let's introduce the following notations:\\

{\bf Notations.} 

1) Denote by $Z_i=span\{x_{i-1}, \ x_{i}\}, \ i\in\mathbb{Z}_n.$ \\
Clearly, $Z_i\subseteq Im(A_i);$  $Z_i\subseteq U,$ and by Lemma 5.1.6(1),
$dim\ Z_i=2.$\\

2) For $i\neq j, i\neq j+1,$  $i, j \in \mathbb{Z}_n,$  denote by\\ $Y_{i,j}=span \{v_{i,j}, v_{i,j+1}\}.$\\

3) For $i\in\mathbb{Z}_n,$ denote by $X^{\uparrow}_{i}=Y_{i,i+1}=span \{x_{i}, v_{i,i+2}\};$ and\\
 $X^{\downarrow}_i=Y_{i,i-2}=span \{x_{i-1}, v_{i,i-2}\};$ \\

{\bf Lemma 5.1.16.} Under conditions {\bf ($\bigstar$)}  for $n\geqslant 5,$ for every $i\in \mathbb{Z}_n,$ \\ $v_{i-1,i+1}\notin Im(A_i).$

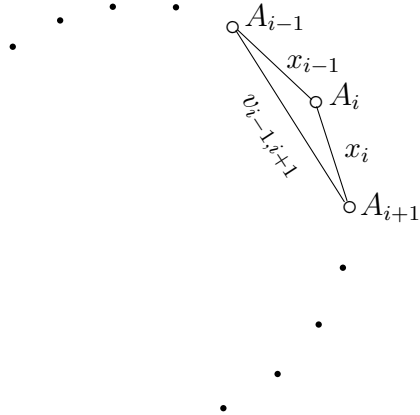
\setlength{\unitlength}{0.7cm}
\begin{picture}(10,10)(-5,-5)
\put(3.85,0){$\circ$}
\put(4.2,0){$A_{i+1}$}
\put(3.21,2){$\circ$}
\put(3.6,2.1){$A_{i}$}
\put(1.63,3.42){$\circ$}
\put(2,3.6){$A_{i-1}$}

\put(3.9,0.25){\line(-1,1.58){2.065}}

\put(3.95,0.27){\line(-1,3.1){0.58}}

\put(1.87,3.52){\line(1,-0.94){1.41}}

\put(0.7,3.95){\circle*{0.12}}
\put(-0.5,3.96){\circle*{0.12}}
\put(-1.5,3.7){\circle*{0.12}}
\put(-2.4,3.2){\circle*{0.12}}

\put(2.8,2.8){\small{$x_{i-1}$}}
\put(3.9,1.1){\small{$x_{i}$}}
\put(1.8,2.2){{\rotatebox{303}{\small{$v_{i-1,i+1}$}}}}

\put(3.87,-1){\circle*{0.12}}
\put(3.41,-2.08){\circle*{0.12}}
\put(2.63,-3.02){\circle*{0.12}}
\put(1.6,-3.67){\circle*{0.12}}

%\put(0,0){\circle{8}}
%\put(0,0){\circle*{0.08}}

\end{picture}

\begin{proof}

 Due to the cyclic argument, it is enough to show that \\$v_{1,3}\notin Im(A_2).$\\

Consider \\$Im(A_1)\cap Im(A_2)\cap Im(A_3)=\big(Im(A_1)\cap Im(A_2)\big)\cap \big( Im(A_3)\cap Im(A_2)\big)=$\\$=span\{x_1\}\cap span\{x_2\} =\{0\},$ since $x_1$ and $x_{2}$ are linearly independent  by Lemma 5.1.6, part 1. Thus,\\

$\{0\}=Im(A_1)\cap Im(A_2)\cap Im(A_3)=\big(Im(A_1)\cap Im(A_3)\big)\cap Im(A_2)=$\\
$=span\{v_{1,3}\}\cap Im(A_2),$ so $v_{1,3}\notin Im(A_2).$

 \end{proof}
  
{\bf Corollary 5.1.17.}  Under conditions {\bf ($\bigstar$)} for $n\geqslant 5,$  \\for every $i\in \mathbb{Z}_n,$\\
1) $v_{i-1,i+1}\notin Z_i=span\{x_{i-1}, \ x_{i}\};$\\
2) $x_{i}$ and  $v_{i,i+2}$ are linearly independent; \\
3) $dim(X^{\uparrow}_{i})=2;$\\
4) $x_{i-1}$ and  $v_{i-2,i}$ are linearly independent; \\
5) $dim(X^{\downarrow}_{i})=2;$\\

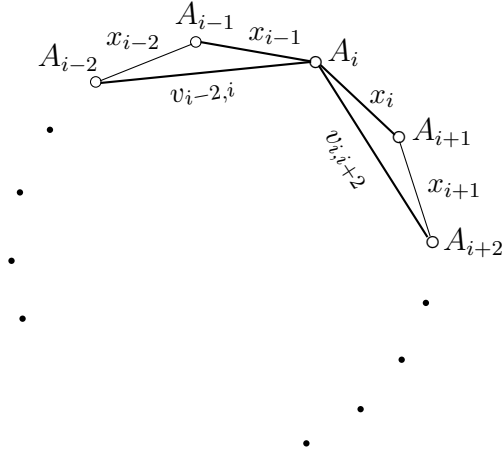
\setlength{\unitlength}{0.7cm}
\begin{picture}(10,10)(-5,-5)
\put(3.85,0){$\circ$}
\put(4.2,0){$A_{i+2}$}
\put(3.21,2){$\circ$}
\put(3.6,2.1){$A_{i+1}$}
\put(1.63,3.42){$\circ$}
\put(2,3.6){$A_{i}$}

\thicklines
\put(3.9,0.25){\line(-1,1.58){2.065}}
\thinlines
\put(3.95,0.27){\line(-1,3.1){0.58}}
\thicklines
\put(1.87,3.52){\line(1,-0.94){1.41}}
\thinlines

\put(-0.5,3.96){\circle{0.2}}
\put(-0.9,4.3){$A_{i-1}$}
\put(-2.4,3.2){\circle{0.2}}
\put(-3.5,3.5){$A_{i-2}$}

\put(2.8,2.8){\small{$x_{i}$}}
\put(3.9,1.1){\small{$x_{i+1}$}}
\put(1.8,2.2){{\rotatebox{303}{\small{$v_{i,i+2}$}}}}

\put(3.87,-1){\circle*{0.12}}
\put(3.41,-2.08){\circle*{0.12}}
\put(2.63,-3.02){\circle*{0.12}}
\put(1.6,-3.67){\circle*{0.12}}

\thicklines
\put(1.69,3.59){\line(-1,0.175){2.1}}
\put(1.69,3.58){\line(-1,-0.10){4.0}}
\thinlines
\put(-2.32,3.22){\line(1,0.40){1.75}}

\put(-1,2.8){{\rotatebox{10}{\small{$v_{i-2,i}$}}}}
\put(0.5,4){\small{$x_{i-1}$}}
\put(-2.2,3.9){\small{$x_{i-2}$}}

\put(-3.27,2.29){\circle*{0.12}}
\put(-3.84,1.1){\circle*{0.12}}
\put(-4,-0.2){\circle*{0.12}}
\put(-3.79,-1.3){\circle*{0.12}}

%\put(0,0){\circle{8}}
%\put(0,0){\circle*{0.08}}

\end{picture}

{\bf Lemma 5.1.18.} Under conditions {\bf ($\bigstar$)} for $n\geqslant 7,$ the vectors $v_{i,j}$ and $v_{i,j+1}$ are linearly independent  for all  $i\neq j, i\neq j+1,$  (indices taken modulo $n$).

\begin{proof}

 Due to the cyclic argument, it is enough to prove  linear independence of vectors $v_{1,k}$ and $v_{1,k+1}$ for all $k=2, \ 3, \ \dots , n-1.$\\

For $k=2$ and $k=n-1,$ the statement is equivalent to  parts 2) and 4) of  Corollary 5.1.17.\\

Suppose now that $3\leqslant k \leqslant n-2,$ and the vectors $v_{1,k}$ and $v_{1,k+1}$ are linearly dependent. The conjugation by $D=\rho(\Delta)$ followed by the conjugation by $T^{k+2}$ sends the vertices $A_1, \ A_{k}$ and $A_{k+1}$ to $A_{k+1}, \ A_{2}$ and $A_1$ respectively.

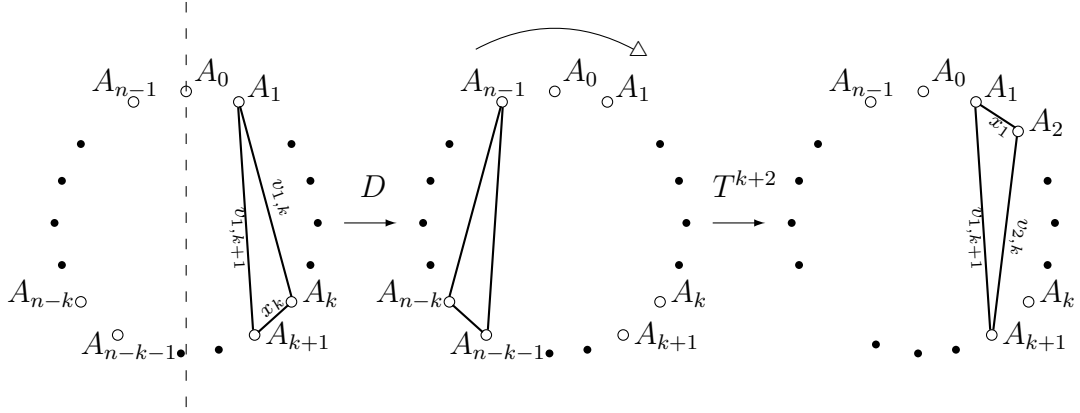
\setlength{\unitlength}{0.7cm}
\begin{picture}(20,10)(-10,-5)

\multiput(-5,0)(7,0){2}{\vector(1,0){1}}

\multiput(-8,-3.5)(0,0.5){16}{\line(0,0.3){0.2}}

\qbezier(-2.5,3.3)(-1,4.1)(.5,3.3)
\put(.4,3.17){\tiny{$\triangle$}}
\put(-4.7,0.5){$D$}
\put(2.,0.5){$T^{k+2}$}

\put(7.8,1.74){\circle{0.2}}
\put(8,1.74){$A_2$}

\multiput(-8,2.5)(7,0){3}{\circle{0.2}}
\multiput(-7.85,2.7)(7,0){3}{$A_0$}
\multiput(-7,2.3)(7,0){3}{\circle{0.2}}
\multiput(-6.85,2.4)(7,0){3}{$A_1$}
\multiput(-9,2.3)(7,0){3}{\circle{0.2}}
\multiput(-9.8,2.5)(7,0){3}{$A_{n-1}$}
\multiput(-6,-1.5)(7,0){3}{\circle{0.2}}
\multiput(-5.8,-1.5)(7,0){3}{$A_{k}$}
\multiput(-10,-1.5)(7,0){2}{\circle{0.2}}
\multiput(-11.4,-1.5)(7,0){2}{$A_{n-k}$}
\multiput(-6.7,-2.13)(7,0){3}{\circle{0.2}}
\multiput(-6.5,-2.3)(7,0){3}{$A_{k+1}$}
\multiput(-9.3,-2.13)(7,0){2}{\circle{0.2}}
\multiput(-10,-2.5)(7,0){2}{$A_{n-k-1}$}

\thicklines
\put(-6,-1.4){\line(-1,3.7){0.978}}
\put(-6.06,-1.56){\line(-1,-0.9){0.575}}
\multiput(-6.73,-2.03)(14,0){2}{\line(-1,15){0.282}}

\put(-3,-1.4){\line(1,3.7){0.978}}
\put(-2.27,-2.03){\line(1,15){0.282}}
\put(-2.94,-1.56){\line(1,-0.9){0.575}}
\put(7.33,-2.03){\line(1,8.2){.448}}
\put(7.75,1.81){\line(-1,.65){.67}}

\thinlines

\put(-6.5,0.8){{\rotatebox{290}{\tiny{$v_{1,k}$}}}}
\put(7.6,0.1){{\rotatebox{264}{\tiny{$v_{2,k}$}}}}
\put(-7.2,0.3){{\rotatebox{273}{\tiny{$v_{1,k+1}$}}}}
\put(6.8,0.3){{\rotatebox{273}{\tiny{$v_{1,k+1}$}}}}
\put(-6.6,-1.8){{\rotatebox{45}{\tiny{$x_{k}$}}}}
\put(7.2,1.9){{\rotatebox{325}{\tiny{$x_{1}$}}}}

\multiput(-5.5,0)(7,0){3}{\circle*{0.15}}
\multiput(-5.64,0.8)(7,0){3}{\circle*{0.15}}
\multiput(-6.,1.5)(7,0){2}{\circle*{0.15}}
\multiput(-5.64,-0.8)(7,0){3}{\circle*{0.15}}

\multiput(-10.5,0)(7,0){3}{\circle*{0.15}}
\multiput(-10.36,0.8)(7,0){3}{\circle*{0.15}}
\multiput(-10,1.5)(7,0){3}{\circle*{0.15}}
\multiput(-10.36,-0.8)(7,0){3}{\circle*{0.15}}

\multiput(-8.1,-2.48)(7,0){3}{\circle*{0.15}}
\multiput(-7.38,-2.41)(7,0){3}{\circle*{0.15}}
\put(5.1,-2.31){\circle*{0.15}}

%\multiput(-8,0)(7,0){3}{\circle{5}}

\end{picture}

Thus, if  $span\{v_{1,k}\}=span\{v_{1,k+1}\}$ then \\ $span\{v_{1,k+1}\}=span\{v_{2,k+1}\}.$ \\By Lemma 4.16, we obtain \\
$span\{v_{1,k+1}\}=span\{x_k\}=span\{x_{1}\},$ a contradiction with Corollary 5.1.13, since $x_1$ and $x_k$ are linearly independent for $n\geqslant 7.$

\end{proof}

{\bf Corollary 5.1.19.} Under conditions {\bf ($\bigstar$)} for $n\geqslant 7,$\\$dim(Y_{i,j})=2$ for all $i\neq j, i\neq j+1,\ i,j\in\mathbb{Z}_n.$\\

{\bf Lemma 5.1.20.} Let $U_1=Im(A_1)$ and $U_k=U_{k-1}+Im(A_k)$ for \\$k=2, \ 3, \ \dots , n-1.$\\

If $n\geqslant 5,$ and the conditions {\bf ($\bigstar$)} hold, then $dim(U_k)\leqslant k+3.$

\begin{proof}

$dim(U_1)=3, $ and \\$dim(U_2)=dim(U_1)+dim(Im(A_2))-dim(U_1\cap Im(A_2))=$\\$=dim(Im(A_1))+dim(Im(A_2))-dim\big(Im(A_1)\cap Im(A_2)\big)=$\\$=3+3-1=5.$\\

By definition, $X^{\downarrow}_{k}\subseteq Im(A_k)$ and $X^{\downarrow}_{k}\subseteq Im(A_{k-1})+Im(A_{k-2}),$ so \\$X^{\downarrow}_{k}\subseteq U_{k-1}\cap Im(A_k)$ for $k=3, \ 4, \ \dots, n-1.$ By Corollary 5.1.17, part 5), since $n\geqslant 5,$ $dim(X^{\downarrow}_{k})=2,$ and thus, $dim\big(U_{k-1}\cap Im(A_k)\big)\geqslant 2.$\\

By induction on $k$ we have:\\ 
$dim(U_k)=dim(U_{k-1})+dim(Im(A_k))-dim(U_{k-1}\cap Im(A_k))\leqslant$\\$\leqslant [(k-1)+3]+3-2=k+3.$
\end{proof}

	{\bf Lemma 5.1.21.}  Under conditions {\bf ($\bigstar$)} for $n\geqslant 5,$ for all $i\in \mathbb{Z}_n$\\
1)  $v_{i,k}\in X^{\uparrow}_{i} $ for $k\in \mathbb{Z}_n,$ $k\neq i-2, \ i-1,\  i; $  \\
2)  $v_{i,k}\in X^{\downarrow}_{i} $ for $k\in \mathbb{Z}_n,$ $k\neq i, \ i+1,\  i+2.$

\begin{proof}

  Due to cyclicity, to prove part 1), it is enough to prove that  $v_{1,k}\in X^{\uparrow}_{1} $ for $k\neq n-1, \ 0,\ 1. $ Part 2) follows from part 1) by applying the conjugation by $D$ followed by the conjugation by $T^2.$\\

For $k=2$ and $k=3$ the statement $v_{1,k}\in X^{\uparrow}_{1} $ follows from the definition of $X^{\uparrow}_{1}.$  \\

Suppose now for some  $k, \ k=4, \ 5, \ \dots , n-2,$ \\$v_{1,k}\notin X^{\uparrow}_{1}. $ Then $dim(span\{v_{1,k}, \ x_1, \ v_{1,3}\})=3,$ and\\
  $Im(A_1)=span\{v_{1,k}, \ x_1, \ v_{1,3}\}.$\\ Conjugating by $T^{k+1}D$ gives \\
 $Im(A_k)=span\{v_{1,k}, \ x_{k-1}, \ v_{k,k-2}\}.$

\setlength{\unitlength}{0.7cm}
\begin{picture}(20,10)(-10,-5)

\multiput(-5,0)(7,0){2}{\vector(1,0){1}}

\multiput(-8,-3.5)(0,0.5){16}{\line(0,0.3){0.2}}

\qbezier(-2.5,3.3)(-1,4.1)(.5,3.3)
\put(.4,3.17){\tiny{$\triangle$}}
\put(-4.7,-0.6){$D$}
\put(2.,0.2){$T^{k+1}$}

\multiput(-6.15,1.7)(7,0){3}{\circle{0.2}}
\multiput(-6,1.74)(7,0){3}{$A_2$}
\put(-5.64,0.8){\circle{0.2}}
\put(-5.5,0.8){$A_3$}

\put(-2.85,1.7){\circle{0.2}}
\put(-4,1.9){$A_{n-2}$}
\put(-3.36,0.8){\circle{0.2}}
\put(-4.4,1.05){$A_{n-3}$}

\multiput(-8,2.5)(7,0){3}{\circle{0.2}}
\multiput(-7.85,2.7)(7,0){3}{$A_0$}
\multiput(-7,2.3)(7,0){3}{\circle{0.2}}
\multiput(-6.85,2.4)(7,0){3}{$A_1$}
\multiput(-9,2.3)(7,0){3}{\circle{0.2}}
\multiput(-9.8,2.5)(7,0){3}{$A_{n-1}$}
\put(1,-1.5){\circle{0.2}}
\put(1.2,-1.5){$A_{n-k}$}
\multiput(-10,-1.5)(14,0){2}{\circle{0.2}}
\multiput(-10.8,-1.7)(14,0){2}{$A_{k}$}

\put(-6.62,2.17){{\rotatebox{330}{\tiny{$x_{1}$}}}}
\put(-6.9,1.8){{\rotatebox{312}{\tiny{$v_{1,3}$}}}}
\multiput(-8.9,0.3)(14,0){2}{{\rotatebox{58}{\tiny{$v_{1,k}$}}}}

\thicklines
\put(-6.2,1.8){\line(-1,.67){.71}}
\put(-6.96,2.23){\line(1,-1.08){1.26}}
\multiput(-7.08,2.23)(14,0){2}{\line(-1,-1.28){2.86}}
\put(-2.75,1.75){\line(1,.72){.68}}
\put(-2.04,2.23){\line(-1,-1.08){1.26}}
\put(-1.92,2.23){\line(1,-1.28){2.86}}
\put(4.05,-1.6){\line(1,-0.79){.59}}
\put(4.1,-1.5){\line(1,-0.6){1.51}}

\thinlines

\put(3.85,-1.86){{\rotatebox{325}{\tiny{$x_{k-1}$}}}}
\put(4.4,-1.5){{\rotatebox{332}{\tiny{$v_{k,k-2}$}}}}

\multiput(-5.5,0)(7,0){3}{\circle*{0.15}}
\multiput(-6.7,-2.13)(7,0){3}{\circle*{0.15}}
\multiput(-6.,-1.5)(14,0){2}{\circle*{0.15}}
\multiput(-5.64,-0.8)(7,0){3}{\circle*{0.15}}
\multiput(1.36,0.8)(7,0){2}{\circle*{0.15}}
\multiput(-10.5,0)(7,0){3}{\circle*{0.15}}
\multiput(-10.36,0.8)(14,0){2}{\circle*{0.15}}
\multiput(-10,1.5)(14,0){2}{\circle*{0.15}}
\multiput(-10.36,-0.8)(7,0){3}{\circle*{0.15}}
\multiput(-9.3,-2.13)(7,0){2}{\circle*{0.15}}
\multiput(-8.4,-2.47)(7,0){2}{\circle*{0.15}}
\multiput(-7.5,-2.45)(7,0){3}{\circle*{0.15}}

\put(-2.96,-1.55){\circle*{0.15}}

\put(4.7,-2.13){\circle{0.2}}
\put(3.5,-2.7){$A_{k-1}$}
\put(5.7,-2.47){\circle{0.20}}
\put(5.2,-3.){$A_{k-2}$}

%\multiput(-8,0)(7,0){3}{\circle{5}}

\end{picture}

 Consider subspaces $U_1, \ U_2, \ \dots , U_{n-1}$ defined by $U_1=Im(A_1)$ and \\ $U_i=U_{i-1}+Im(A_i)$ for $i=2, \dots, n-1.$ \\By Lemma 5.1.20, we have
$dim(U_{k-1})\leqslant k+2.$\\ 
 
 But $v_{1,k}\in Im(A_{1}), \ x_{k-1}\in Im(A_{k-1}),$ and $ v_{k,k-2}\in Im(A_{k-2}),$ so\\ $Im(A_k)\subseteq U_{k-1},$ and \\$dim(U_k)=dim(U_{k-1})\leqslant k+2.$\\
 
 Moreover,  for every $j=k, \ k+1, \dots , n-1,$ by using conjugation by $T^{j-k},$we obtain that  $Im(A_j)=span\{v_{j,j-k+1}, \ x_{j-1}, \ v_{j,j-2}\}\subseteq U_{j-1},$ (since $j-k+1\leqslant j-1$ for $k\geqslant 2$),  and by induction on $j$ we get \\$U_j= U_{j-1}=\dots =U_{k-1},$  thus  $dim(U_j)\leqslant k+2$ for all \\$j=k, \ k+1, \dots , n-1.$\\
 
 On the other hand, $U_{n-1}=Im(A_1)+Im(A_2)+\dots+Im(A_{n-1})=V,$ so for $k\leqslant n-2, $ we have \\$n+1\leqslant dim(V)=dim(U_{n-1})\leqslant k+2\leqslant n-2+2=n,$ a contradiction.
 
\end{proof}

 {\bf Lemma 5.1.22.}  Under conditions {\bf ($\bigstar$)} for $n\geqslant 7,$ for all $i\in \mathbb{Z}_n,$\\
 $Z_i=X^{\uparrow}_{i}=X^{\downarrow}_{i}.$ 
 
 \begin{proof} 
 
 As before, it is enough to prove this statement for $i=1$ due to cyclicity.\\
 
 Consider $Y_{1,4}=span \{v_{1,4}, v_{1,5}\}.$ By Corollary 5.1.19, since $n\geqslant 7,$\\ $dim(Y_{1,4})=2.$\\
 By Lemma 5.1.21 part 1), $Y_{1,4}\subseteq X^{\uparrow}_{1},$ since for $k=4$ and $5,$  $k\leq n-2$ for $n\geqslant 7.$ \\ By Corollary 5.1.17, part 3), $dim(X^{\uparrow}_{1})=2,$ so  $Y_{1,4}= X^{\uparrow}_{1}.$\\
 Similarly, by Corollary 5.1.19, Lemma 5.1.21, part 2) and Corollary 5.1.17,  part 5), $Y_{1,4}= X^{\downarrow}_{1}.$\\
 Thus, $X^{\uparrow}_{1}= X^{\downarrow}_{1}=Y_{1,4}.$\\
 
 Now, $x_1\in X^{\uparrow}_{1}=Y_{1,4}$ and $x_0\in X^{\downarrow}_{1}=Y_{1,4}.$ Thus,\\ $Z_1=span\{x_0, x_1\}\subseteq Y_{1,4}.$ But \\$dim(Z_1)=2$ and $dim(Y_{1,4})=2,$ so\\
 $Z_1=Y_{1,4}=X^{\uparrow}_{1}= X^{\downarrow}_{1}.$
  \end{proof}

 {\bf Corollary 5.1.23.}  Under conditions {\bf ($\bigstar$)} for $n\geqslant 7$, $v_{i,k}\in Z_i$ for all $i,k\in \mathbb{Z}_n,$ $ k\neq i.$

  \begin{proof}
  
   For all $k\neq  i-2, \ i-1, \  i+1, \ i+2, $ by Lemma 5.1.21, \\$v_{i,k}\in X^{\uparrow}_{1}\cap X^{\downarrow}_{1}=Z_i$\\
  For $k=i+1, \ i+2, $ $v_{i,k}\in X^{\uparrow}_{1}=Z_i$ and for $k=i-1, \ i-2,$ $ v_{i,k}\in X^{\downarrow}_{1}=Z_i$ by Lemma 5.1.22.
  
   \end{proof}

  {\bf Theorem 5.1.24.} Let $\rho: B_n\to GL_{r}(\mathbb{C})$ be an {\it irreducible} matrix representation of $B_n$ of dimension $r\geqslant n+1$  with $dim(Im(A_i))=3$ for $n\geqslant 9.$\\
  
  Then $Im(A_i)\cap Im(A_{i+1})=\{0\}$ for all $i\in \mathbb{Z}_n.$
  
   \begin{proof}
   
    Suppose $Im(A_i)\cap Im(A_{i+1})\neq\{0\}.$ Then by Theorem 5.1.15,
   $dim(Im(A_i)\cap Im(A_{j}))= 1$  for all $i, j\in \mathbb{Z}_n, \ i\neq j.$\\
   
   Consider $v_{1,3}=v_{3,1}\neq 0,$ such that $span\{v_{1,3}\}=Im(A_1)\cap Im(A_{3}).$
   By Corollary 5.1.23, $v_{1,3}\in Z_1=span\{x_0, x_1\},$ and  \\$v_{3,1}\in Z_3=span\{x_2, x_3\}.$
  So, $v_{1,3}=v_{3,1}\in span\{x_0, x_1\}\cap span\{x_2, x_3\}.$\\Since $n\geqslant 9, $ the set of vectors $x_0,\  x_1, \ x_2, \ x_3$ is a subset of a linearly independent set $x_0, \ x_1, \ \dots, x_{n-4}$ (Corollary 5.1.12), thus $ v_{1,3}=0,$ a contradiction.
  
  \end{proof}
  
  { \centerline{5.2. \bf{ Case II: $Im(A_i)\cap Im(A_{i+1})= \{0\}.$}}}
 
  \vskip 0.3cm 

In this subsection we are going to follow the path similar to the case when $Im(A_i)\cap Im(A_{i+1})\neq \{0\}.$  First, we are going to prove that the friendship graph has all non-adjacent edges (diagonals) and that each diagonal represents a one-dimensional subspace (Lemma 5.2.1). Then we are going to show (Theorems 5.2.13 and 5.2.14) that there are no irreducible  representations associated to such friendship graphs if the number of vertices is at least $10.$\\
 
     {\bf Lemma 5.2.1.} Let $\rho: B_n\to GL_{r}(\mathbb{C})$ be an {\it irreducible} matrix representation of $B_n$ of dimension $r\geqslant n+1,$ $n\geqslant 5,$ with \\$dim(Im(A_i))=3.$ Suppose that $Im(A_i)\cap Im(A_{i+1})=\{0\}$ for all $i\in \mathbb{Z}_n.$ \\Then $dim(Im(A_i)\cap Im(A_{j}))=1$ for all $||i-j||\geqslant 2, \  i, \  j \in \mathbb{Z}_n.$
 
 \begin{proof} Due to the cyclicity, it is enough to show that \\$dim(Im(A_1)\cap Im(A_{k}))=1$ for all $k=3, \ 4, \dots, n-1.$\\
 
 By Corollary 4.11, $dim(Im(A_1)\cap Im(A_{k}))\geqslant 1$ 
  and\\
 by  Lemma 4.14,  $dim(Im(A_1)\cap Im(A_{k}))\leqslant 2$ for $k=3, \ 4, \dots, n-1.$ 
 \\
 By Lemma 4.13, $dim(Im(A_1)\cap Im(A_{k}))=dim(Im(A_1)\cap Im(A_{m}))$ \\for all $k, \ m=3, \ 4, \ \dots n-1.$ Thus, it is enough to prove that  \\$dim(Im(A_1)\cap Im(A_{k}))\neq 2.$ \\
 
 Suppose $dim(Im(A_1)\cap Im(A_{k}))= 2.$ Then, \\ $dim(Im(A_2)\cap Im(A_{4}))= 2,$ and \\
 
  $Im(A_1)\cap Im(A_{4})+Im(A_2)\cap Im(A_{4}) \subseteq Im(A_4),$ so \\

$dim( Im(A_1)\cap Im(A_{2}))\geqslant dim(Im(A_1)\cap Im(A_{2})\cap Im(A_{4}))=$

$=dim(Im(A_1)\cap Im(A_{4}))+dim(Im(A_2)\cap Im(A_{4}))-$

 $-dim(Im(A_1)\cap Im(A_{4})+Im(A_2)\cap Im(A_{4}))\geqslant$
 
 $\geqslant  2+2-dim( Im(A_{4}))=1,$\\
 
a contradiction with  $Im(A_1)\cap Im(A_{2})=\{0\}.$\\

Thus, $dim(Im(A_1)\cap Im(A_{k}))=1$ for all $k=3, \ 4, \dots, n-1.$

\end{proof}

 To investigate the friendship graphs with one-dimensional diagonals, we are going to use the following notations:\\

{\bf Notations.}

1) For $||i-j||\geqslant 3,$ denote by\\ $X_{i,j}=span \{v_{i,j-1}, \ v_{i,j}, \ v_{i,j+1}\}.$\\

2) Similar to the notation used  in subsection 5.1, denote by \\$Y_{i,j}=span \{ v_{i,j}, \ v_{i,j+1}\}$ for   $i, j \in \mathbb{Z}_n,$ $ j\neq i-2, \ i-1, \  i,\  i+1.$\\

3) Denote by $V_i$  a subspace spanned by all diagonals coming out of $A_i:$ \\$V_i=span \{v_{i,i+2}, \ v_{i,i+3}, \ \dots ,v_{i,i-2}  \}.$\\

 4) Denote by $W$ the subspace of $V$ generated by all the diagonals: \\
 $W=V_1+V_2+\dots +V_{n-1}+V_0.$\\

   {\bf Lemma 5.2.2.} Let $\rho: B_n\to GL_{r}(\mathbb{C})$ be a matrix representation \\of $B_n$ for $n\geqslant 5.$ Suppose that $Im(A_i)\cap Im(A_{i+1})=\{0\}$ for all $i\in \mathbb{Z}_n,$ and $dim(Im(A_i)\cap Im(A_{j}))=1$ for all $||i-j||\geqslant 2, \  i, \  j \in \mathbb{Z}_n.$\\ 
   
 Then for all $i, \ j\in \mathbb{Z}_n, \ j\neq i-2, \ i-1, \  i,\  i+1,$  the vectors $v_{i,j}$ and $v_{i,j+1}$ are linearly independent.

\begin{proof}

  Suppose that $v_{i,j}$ and $v_{i,j+1}$ are linearly dependent for some $i, j.$ 
 Then $v_{i,j+1}\in span\{v_{i,j}\}=Im(A_i)\cap Im(A_{j}),$ so \\ $0\neq v_{i,j+1}\in  Im(A_{j})\cap Im(A_{j+1}),$ a contradiction with \\$Im(A_j)\cap Im(A_{j+1})=\{0\}.$
 
 \end{proof}
 
    {\bf Corollary 5.2.3.} Under conditions of Lemma 5.2.2, \\$dim X_{i,j} \geqslant 2$ for all $||i-j||\geqslant 3.$\\
    
     {\bf Corollary 5.2.4.} Under conditions of Lemma 5.2.2,\\$dim Y_{i,j}= 2$ for all $i,\ j\in \mathbb{Z}_n, \ j\neq i-2, \ i-1, \ i, \ i+1.$\\

    {\bf Lemma 5.2.5.} Suppose $n\geqslant 7.$  Then there exists $k$ such that \\
$1<k<\frac{n}{2}$ and $gcd(k,n)=1$.

\begin{proof}

 Suppose for every $k,\  1<k<\frac{n}{2},\  $ $gcd(k,n)\neq 1$ . \\
Since $gcd(n-k,n)=gcd(k,n)$ and $\frac{n}{2}<n-k<n-1,$  the only integers $k$ in the range $1\leqslant k \leqslant n $ coprime with $n$  are $1$ and $n-1.$ \\
In other words, if $\varphi$ is the
Euler's function, then $\varphi(n)=2.$

 If $n=\prod p_i^{k_i}$ is the prime factorization of
$n$, then by the well-known formula for the Euler function $\varphi(n)=\prod
(p_i-1)p_i^{k_i-1}$. \\Since $\varphi(n)=2,$ the only prime
 factors of $n$ are $2$ and $3$. Moreover, the power of $3$ is at most
$1$ and the power of $2$ is at most $2.$ Thus, $n$ divides $12.$

Since $n\geqslant 7,$ then $n=12.$ But $\varphi(12)=4,$ a contradiction with  $\varphi(n)=2.$ 

 \end{proof}

 {\bf Lemma 5.2.6.} Let $\rho: B_n\to GL_{r}(\mathbb{C})$ be a matrix representation \\of $B_n$ of dimension $r,$ $n\geqslant 9.$ Suppose that $Im(A_i)\cap Im(A_{i+1})=\{0\}$ for all $i\in \mathbb{Z}_n,$ and $dim(Im(A_i)\cap Im(A_{j}))=1$ for all $||i-j||\geqslant 2, \  i, \  j \in \mathbb{Z}_n.$\\ Suppose that $dim(V_1)=2.$
   \\ Then $W=V_1+V_2+\dots +V_{n-1}+V_0$ is $B_n-$invariant.
   
   \begin{proof} 
   
   Due to the cyclic argument, because $V_1$ is 2-dimensional,  $V_i$ is 2-dimensional for all $i\in \mathbb{Z}_n.$\\
   
    For every fixed  $i\in \mathbb{Z}_n,$ by Lemma 5.2.2, the vectors $v_{i,j}$ and $v_{i,j+1}$ are linearly independent for every  $ j\in \mathbb{Z}_n, \ j\neq i-2, \ i-1, \ i, \ i+1,$ and since $dim (V_i)=2,$ we have that $V_i=span\{v_{i,j}, \ v_{i,j+1}\}$ for $ i, \ j\in \mathbb{Z}_n,$\\$ j\neq i-2, \ i-1, \ i, \ i+1.$\\

   To show that $W$ is $B_n-$ invariant, we need to show that for every diagonal $v_{i,j},$ and  for every $k\in \mathbb{Z}_n,$ $A_k v_{i,j}\in W.$  Due to cyclicity, it is enough to show that $A_1 v_{i,j}\in W$ for all $v_{i,j}.$ Depending on the values of $i, \ j,$ we have to consider the following cases:\\

   a) $v_{i,j}$ is coming out of $A_1.$\\
   In this case, $v_{i,j}=v_{1,j}\in Im(A_1)\cap Im(A_{j}), $ where $j=3,\ 4, \dots , n-1.$ Since $A_1$ and $A_j$ commute for $j=3,\ 4, \dots , n-1,$ we have \\$A_1v_{1,j}\in Im(A_1)\cap Im(A_{j})=span\{ v_{1,j}\} \subseteq W.$\\

   b) $v_{i,j}$ connects the vertices $A_i$ and $A_j,$ where both $A_i$ and $A_j$ are not neighbors of $A_1,$ that is $i,\ j=3,\ 4,\ \dots n-1, \ ||i-j||\geqslant 2.$ \\ In this case, $v_{i,j}\in Im(A_i)\cap Im(A_{j}), $ and since $A_1$ commutes with both $A_i$ and $A_j,$  we have $A_1v_{i,j}\in Im(A_i)\cap Im(A_{j})=span\{ v_{i,j}\} \subseteq W.$\\

   c) $v_{i,j}$ connects the vertices $A_i$ and $A_j,$ where $A_i$ is a neighbor of $A_1$ and $A_j$ is not a neighbor of $A_1.$ So, $i=2,\  j=4,\ 5, \dots, \ n-1,$ or\\ $i=0,\  j=3, \ 4, \dots, \ n-2.$\\ Since $n\geqslant 9,$ we can pick 2 vertices $A_k$ and $A_{k+1},$ such that they both are not neighbors of $A_1$ and not neighbors of $A_j.$ (For example, for $j\leqslant 5,$ we can pick the vertices $A_7$ and $A_8,$ and for $j\geqslant 6,$ we can pick $A_3$ and $A_4.$)\\

  Consider  the diagonals $v_{j,k}$ and $v_{j,k+1}.$  Since they form a basis of $V_j,$ \\$v_{i,j}\in span\{v_{j,k}, \ v_{j,k+1}\}.$ Since $A_j,$ $A_k$ and $A_{k+1}$ are not neighbors of $A_1,$ by case (b) we have that
   $A_1v_{j,k}\in W$ and $A_1v_{j,k+1}\in W,$ and hence  $A_1v_{i,j}\in W.$ \\

   d) $v_{i,j}$ connects the vertices $A_{0}$ and $A_{2},$ that is $v_{i,j}=v_{0,2}.$\\
   Since $n\geqslant 9, $ then $A_{3}$ and $A_{4}$ are not neighbors of $A_{0},$ (and not neighbors of $A_1$). The diagonals $v_{0,3}$ and $v_{0,4}$ form a basis of $V_{0},$  so $v_{0,2}\in span\{v_{0,3}, \ v_{0,4}\},$ thus $A_1v_{0,2}\in W,$ since $A_1v_{0,3}\in W$ and $A_1v_{0,4}\in W$ by case (c).
   
  \end{proof}

      {\bf Lemma 5.2.7.} Let $\rho: B_n\to GL_{r}(\mathbb{C})$ be an {\it irreducible} matrix \\ representation of $B_n$ of dimension $r\geqslant n+1,$  $n\geqslant 9.$ Suppose that  \\$Im(A_i)\cap Im(A_{i+1})=\{0\}$ for all $i\in \mathbb{Z}_n,$ and $dim(Im(A_i)\cap Im(A_{j}))=1$ for all $||i-j||\geqslant 2, \  i, \  j \in \mathbb{Z}_n.$\\ Then  $dim(V_i)\geqslant 3$ for all $i\in \mathbb{Z}_n.$

  \begin{proof} Suppose that $dim(V_i)\leqslant 2$ for some $i\in \mathbb{Z}_n.$ By Lemma 5.2.2, $dim(V_i)\geqslant 2,$ hence, $dim(V_i)=2.$ So, due to cyclicity, $dim(V_i)=2$ for all $i\in \mathbb{Z}_n.$
  Then by Lemma  5.2.6, $W$ is $B_n-$invariant. Since $\rho$ is irreducible, $W=V$ and $dim ( W)=r\geqslant n+1.$\\
   
  Since $n\geqslant 9, $  then by Lemma 5.2.5, there exists a number $k,$\\ $1< k< \frac{n}{2}, $ such that $k$ and $n$ are relatively prime.\\

  We claim that the two diagonals $v_{0,k}$ and $v_{0,n-k}$  form a basis of $V_0,$ so $V_0=span\{v_{0,k},\ v_{0,n-k}\}.$\\
  Suppose not. Then $v_{0,n-k}\in span \{v_{0,k}\},$   so \\ $Im(A_0)\cap Im(A_{k})=Im(A_0)\cap Im(A_{n-k}).$\\
  Consider the sequence of vertices $A_0, \ A_{k}, \ A_{2k}, \dots, \ A_{sk}, \dots ,$ \\where $s\in \mathbb{N}.$\\
  
 By the cyclic argument, for every vertex $A_{sk}$  ($  s\geqslant 2$) in the sequence  we have that  $Im(A_{sk})\cap Im(A_{(s-1)k})=Im(A_{(s-2)k})\cap Im(A_{(s-1)k}),$
   and by using induction on $s$, we have \\$Im(A_{sk})\cap Im(A_{(s-1)k})=Im(A_{0})\cap Im(A_{k}).$
  
 By applying Lemma 4.16 to the vertices $A_0, \ A_{k}, \ A_{2k},$ we get  \\$Im(A_{0})\cap Im(A_{2k})=Im(A_{0})\cap Im(A_{k}),$ and  by using induction on $s$ and applying  Lemma 4.16 to the vertices $A_0, \ A_{(s-1)k}, \ A_{sk},$  we have  that $Im(A_{0})\cap Im(A_{sk})=Im(A_{0})\cap Im(A_{k})$ for all $s\in \mathbb{N}.$

  \vskip -0.5cm
  
\setlength{\unitlength}{0.65cm}
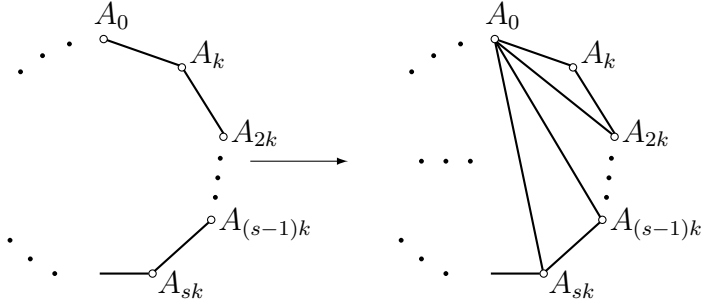
\begin{picture}(20,10)(-10,-5)
\put(-5,0){\vector(1,0){2}}

\multiput(-8,2.5)(8,0){2}{\circle{0.15}}
\multiput(-6.4,1.92)(8,0){2}{\circle{0.15}}
\multiput(-5.55,0.5)(8,0){2}{\circle{0.15}}
\multiput(-5.8,-1.2)(8,0){2}{\circle{0.15}}
\multiput(-7,-2.3)(8,0){2}{\circle{0.15}}

\thicklines
\multiput(-7.08,-2.3)(8,0){2}{\line(-1,0){1}}
\multiput(-6.96,-2.25)(8,0){2}{\line(1,0.9){1.11}}

\multiput(-5.55,0.57)(8,0){2}{\line(-1,1.6){0.811}}
\multiput(-6.45,1.96)(8,0){2}{\line(-1,0.35){1.48}}

\put(0.04,2.44){\line(1,-0.8){2.36}}
\put(0.03,2.43){\line(1,-1.68){2.13}}
\put(0.01,2.43){\line(1,-4.8){0.97}}
\thinlines

\multiput(-0.5,0)(-0.5,0){3}{\circle*{0.1}}

\multiput(-9,-2.3)(8,0){2}{\circle*{0.1}}
\multiput(-9.5,-2)(8,0){2}{\circle*{0.1}}
\multiput(-9.915,-1.618)(8,0){2}{\circle*{0.1}}

\multiput(-8.7,2.4)(8,0){2}{\circle*{0.1}}
\multiput(-9.25,2.16)(8,0){2}{\circle*{0.1}}
\multiput(-9.7,1.83)(8,0){2}{\circle*{0.1}}

\multiput(-8.2,2.8)(8,0){2}{$A_0$}
\multiput(-6.25,2)(8,0){2}{$A_k$}
\multiput(-5.4,0.4)(8,0){2}{$A_{2k}$}
\multiput(-5.65,-1.4)(8,0){2}{$A_{(s-1)k}$}
\multiput(-6.9,-2.7)(8,0){2}{$A_{sk}$}

\multiput(-5.715,-0.74)(0.05,0.4){3}{\circle*{0.1}}
\multiput(2.285,-0.74)(0.05,0.4){3}{\circle*{0.1}}

%\multiput(-5.765,-1.14)(8,0){2}{\line(1,8){0.195}}
%\multiput(-8,0)(8,0){2}{\circle{5}}
\end{picture}
 
 \vskip -0.3cm

  Since $k$ and $n$ are relatively prime, there exists $s\in \mathbb{N},$ such that \\$sk\equiv 1\pmod{n}.$ Thus, $Im(A_{0})\cap Im(A_{1})=Im(A_{0})\cap Im(A_{k})\neq \{0\},$ a contradiction.\\

  Now, for each $i\in \mathbb{Z}_n,$ consider the two diagonals $v_{i,i+k}$ and $v_{i,i-k} $ (indices are taken modulo $n$). Due to the cyclic argument, \\$V_i=span\{v_{i,i+k}, \ v_{i,i-k}\}$ for each $i, $ and hence \\ $W=span\{ v_{0,k}, \ v_{1,k+1}, \dots , v_{n-1,k-1}, v_{0,n-k}, \ v_{1,n-k+1}, \ \dots , v_{n-1,n-k-1}\}.$ \\
  
  Since $v_{i,i+k}=v_{i+k,i},$ each diagonal in the set generating $W$ appears exactly twice, thus  the total number of the diagonals spanning $W$ is $\frac{2\cdot n}{2}=n.$\\

 So, $dim (W)\leqslant n,$ a contradiction with $dim ( W)=r\geqslant n+1.$
 
 \end{proof}   
     
     Notice that the statements 5.2.2 -- 5.2.7  hold true for the representations of any corank $corank(\rho)\geqslant 3.$ Now we are going restrict our considerations to the representations of corank 3, and  consider the representations satisfying the following conditions:\\
    
\begin{tabular}{cl}
 & Let $\rho: B_n\to GL_{r}(\mathbb{C})$ be an {\it irreducible}  matrix \\ 
 {\bf  ($\bigstar$ $\bigstar$ )}& representation of $B_n$ of dimension $r\geqslant n+1$ with $rk(A_1)=3,$  \\
 &  
where $dim(Im(A_i)\cap Im(A_{i+1}))=0$  for all $i\in \mathbb{Z}_n,$ and \\
 & $dim(Im(A_i)\cap Im(A_{j}))=1$ for  $||i-j||\geqslant 2,$ $i, \ j\in\mathbb{Z}_n.$
\end{tabular}\\  
 
 \vskip 0.3cm
    {\bf Corollary 5.2.8.} Let $\rho: B_n\to GL_{r}(\mathbb{C})$ be a representation satisfying conditions  {\bf  ($\bigstar$ $\bigstar$ )} for $n\geqslant 9.$\\ 
    
    Then for all $i,$ $Im(A_i) =V_i,$ and $V=W=V_1+V_2+\dots +V_{n-1}+V_0.$

  \begin{proof}
  
     By Lemma 5.2.7, $dim(V_i)\geqslant 3$ for every $i\in \mathbb{Z}_n.$ Together with \\$V_i\subseteq Im(A_i) $ and $dim(Im(A_i))=3, $ we have that
   $dim(V_i)=3,$  \\$ V_i=Im(A_i),$ and \\$W=V_1+V_2+\dots +V_{n-1}+V_0=$\\ $=Im(A_1)+Im(A_2)+\dots +Im(A_{n-1})+Im(A_0)$ \\  is an invariant subspace of $V,$ and, since $\rho$ is irreducible, $V=W.$
   
 \end{proof}
    
   {\bf Lemma 5.2.9.} Under conditions  ($\bigstar$ $\bigstar$) for $n\geqslant 9,$ for every $i\in \mathbb{Z}_n$ there exists $j\in \mathbb{Z}_n, \ ||i-j||\geqslant 3,$ such that $X_{i,j}=Im(A_i).$
   
 \begin{proof}
 
 Due to cyclicity, it is enough to prove that there exists \\$j=4, \ 5, \dots , n-2,$ such that $X_{1,j}=Im(A_1).$\\
   
   Suppose that for all $j=4, \ 5, \dots , n-2,$ we have $X_{1,j}\neq Im(A_1).$ 
   Since $X_{1,j}\subseteq Im(A_1)$ and $dim(Im(A_1))=3,$ then 
   $dim(X_{1,j})\leqslant 2.$ By Corollary 5.2.3,  $dim X_{i,j} \geqslant 2$ for all $||i-j||\geqslant 3,$
   so  $dim(X_{1,j})=2$ for $j=4, \ 5, \dots , n-2.$\\
   
   By definition, $X_{1,j}=span \{v_{1,j-1}, \ v_{1,j}, \ v_{1,j+1}\},$ and by Lemma 5.2.2, the vectors $ v_{1,j-1}$ and  $v_{1,j}$ are linearly independent. Thus,  for all $j,$ \\$4\leqslant j\leqslant n-2,$ we have 
$v_{1,j+1}\in   span \{v_{1,j-1}, \ v_{1,j}\}=Y_{1,j-1}.$\\
In particular, $v_{1,5}\in  Y_{1,3}.$ \\By using induction on $j,$ we have $v_{1,j+1}\in   span \{v_{1,j-1}, \ v_{1,j}\}\subseteq Y_{1,3}$ for all $j=4, \ 5, \dots, n-2,$ and, hence,\\
$V_1=span\{v_{1,3}, \ v_{1,4}, \ \dots ,v_{1,n-1}  \}\subseteq Y_{1,3},$ so $dim(V_1)\leqslant 2$, a contradiction with Lemma 5.2.7.

 \end{proof}

 {\bf Lemma 5.2.10.} Under conditions  ($\bigstar$ $\bigstar$) for $n\geqslant 9,$ there exists \\$j\in \mathbb{Z}_n, \ 4\leqslant j \leqslant \frac{n}{2}+1,$ such that $X_{1,j}=Im(A_1).$

\begin{proof}

 By Lemma 5.2.9, there exists \\$j=4, \ 5, \dots , n-2,$ such that $X_{1,j}=Im(A_1).$\\
Suppose that $j>\frac{n}{2}+1.$ Using conjugation by $D,$ we get\\ $X_{n-1,n-j}=Im(A_{n-1}),$ and using conjugation by $T^2,$ we get\\ $X_{1,n-j+2}=Im(A_{1}).$ \\So, there exists $k=n-j+2<n-\left(\frac{n}{2}+1\right)+2= \frac{n}{2}+1,$ such that \\$X_{1,k}=Im(A_1).$

 \end{proof}
 
   {\bf Lemma 5.2.11.} Let $U_1=Im(A_1)$ and $U_k=U_{k-1}+Im(A_k)$ for \\$k=2, \ 3, \ \dots , n-1.$\\

If $n\geqslant 5$ and  the conditions ($\bigstar$ $\bigstar$) hold, then $dim(U_k)\leqslant k+5$ for all $k\geqslant 3.$

\begin{proof}

$dim(U_1)=3, $ and \\$dim(U_2)=dim(U_1)+dim(Im(A_2))-dim(U_1\cap Im(A_2))=$\\$=dim(Im(A_1))+dim(Im(A_2))-dim\big(Im(A_1)\cap Im(A_2)\big)=$\\$=3+3-0=6.$\\

Since $v_{1,3}\in Im(A_1)\cap Im(A_3)\subseteq U_2\cap Im(A_3),$ it follows that \\ $dim(U_2\cap Im(A_3))\geqslant 1,$ so\\
$dim(U_3)=dim(U_2)+dim(Im(A_3))-dim(U_2\cap Im(A_3))\leqslant 6+3-1=8.$\\

$Y_{k,1}=span\{v_{k,1}, v_{k,2}\}\subseteq \big( Im(A_{1})+Im(A_{2}) \big) \cap Im(A_k)\subseteq $\\$\subseteq U_{k-1}\cap Im(A_k)$ for every $k=4, \ 5, \ \dots, n-1.$  \\By Corollary 5.2.4, $dim(Y_{k,1})=2,$ thus, $dim\big(U_{k-1}\cap Im(A_k)\big)\geqslant 2.$\\

By induction on $k$ we have:\\ 
$dim(U_k)=dim(U_{k-1})+dim(Im(A_k))-dim(U_{k-1}\cap Im(A_k))\leqslant$\\$\leqslant [(k-1)+5]+3-2=k+5.$

  \end{proof}
  
   {\bf Lemma 5.2.12.} Under conditions  ($\bigstar$ $\bigstar$) for $n\geqslant 9,$ let $j\in \mathbb{Z}_n,$\\ $4\leqslant j \leqslant \frac{n}{2}+1,$ be a number such that $X_{1,j}=Im(A_1)$\\
   
    Let $U_1=Im(A_1)$ and $U_k=U_{k-1}+Im(A_k)$ for $k=2, \ 3, \ \dots , n-1.$\\
    
    Then $U_k=U_j$ for all $k=j+1, \ j+2, \dots , n-1.$
    
    \begin{proof}
    
     For every $k, \ j+1\leqslant k\leqslant n-1,$ by using conjugation by $D,$ followed by conjugation by $T^{k+1},$ we have:\\
    
     $X_{1,j}=Im(A_1)\Longrightarrow X_{n-1,n-j}=Im(A_{n-1})\Longrightarrow$\\
     
     $\Longrightarrow X_{n-1+k+1,n-j+k+1}= Im(A_{n-1+k+1})\Longrightarrow$\\ 
     
     $\Longrightarrow Im(A_{k})=X_{k,k-j+1}=span \{v_{k,k-j}, \ v_{k,k-j+1}, \ v_{k,k-j+2}\}.$\\

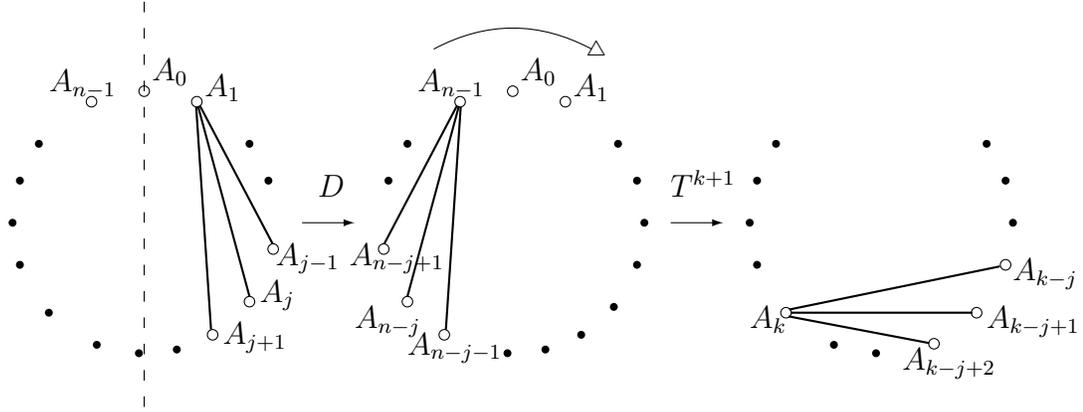
\setlength{\unitlength}{0.7cm}
\begin{picture}(20,10)(-10,-5)

\multiput(-5,0)(7,0){2}{\vector(1,0){1}}

\multiput(-8,-3.5)(0,0.5){16}{\line(0,0.3){0.2}}

\qbezier(-2.5,3.3)(-1,4.1)(.5,3.3)
\put(.4,3.17){\tiny{$\triangle$}}
\put(-4.7,0.5){$D$}
\put(2.,0.5){$T^{k+1}$}

\multiput(-8,2.5)(7,0){2}{\circle{0.2}}
\multiput(-7.85,2.7)(7,0){2}{$A_0$}
\multiput(-7,2.3)(7,0){2}{\circle{0.2}}
\multiput(-6.85,2.4)(7,0){2}{$A_1$}
\multiput(-9,2.3)(7,0){2}{\circle{0.2}}
\multiput(-9.8,2.5)(7,0){2}{$A_{n-1}$}
\put(-6,-1.5){\circle{0.2}}
\put(-5.8,-1.5){$A_{j}$}
\put(-3,-1.5){\circle{0.2}}
\put(-4.,-2){$A_{n-j}$}
\put(-6.7,-2.13){\circle{0.2}}
\put(-6.5,-2.3){$A_{j+1}$}
\put(-2.3,-2.13){\circle{0.2}}
\put(-3,-2.5){$A_{n-j-1}$}
\put(-5.55,-0.5){\circle{0.2}}
\put(-5.5,-0.8){$A_{j-1}$}
\put(-3.45,-0.5){\circle{0.2}}
\put(-4.1,-0.8){$A_{n-j+1}$}

\put(4.19,-1.71){\circle{0.2}}
\put(3.5,-2){$A_k$}
\put(7,-2.3){\circle{0.2}}
\put(6.4,-2.8){$A_{k-j+2}$}
\put(7.81,-1.71){\circle{0.2}}
\put(8,-2){$A_{k-j+1}$}
\put(8.36,-0.8){\circle{0.2}}
\put(8.5,-1.1){$A_{k-j}$}

\thicklines
\put(-6,-1.4){\line(-1,3.7){0.978}}
\put(-6.73,-2.03){\line(-1,15){0.282}}

\put(-3,-1.4){\line(1,3.7){0.978}}
\put(-2.27,-2.03){\line(1,15){0.282}}

\put(-5.56,-0.4){\line(-1,1.9){1.38}}
\put(-3.44,-0.4){\line(1,1.9){1.38}}

\put(4.24,-1.78){\line(1,-0.19){2.67}}
\put(4.29,-1.71){\line(1,0){3.44}}
\put(4.25,-1.65){\line(1,0.205){4}}
\thinlines

\put(0.92,-1.6){\circle*{0.15}}
\put(0.3,-2.13){\circle*{0.15}}
\put(-9.81,-1.71){\circle*{0.15}}

\multiput(1.5,0)(7,0){2}{\circle*{0.15}}
\multiput(-5.64,0.8)(7,0){3}{\circle*{0.15}}
\multiput(-6.,1.5)(7,0){3}{\circle*{0.15}}
\put(1.36,-0.8){\circle*{0.15}}

\multiput(-10.5,0)(14,0){2}{\circle*{0.15}}
\multiput(-10.36,0.8)(7,0){3}{\circle*{0.15}}
\multiput(-10,1.5)(7,0){3}{\circle*{0.15}}
\multiput(-10.36,-0.8)(14,0){2}{\circle*{0.15}}

\multiput(-8.1,-2.48)(7,0){3}{\circle*{0.15}}
\multiput(-7.38,-2.41)(7,0){2}{\circle*{0.15}}
\multiput(5.1,-2.32)(-14,0){2}{\circle*{0.15}}

%\multiput(-8,0)(7,0){3}{\circle{5}}

\end{picture}

\vskip 0.5cm

   We will prove the statement of the lemma by induction on $k.$ \\
   
   First, for $k=j+1,$ we have \\
   $Im(A_{j+1})=span \{v_{j+1,1}, \ v_{j+1,2}, \ v_{j+1,3}\}\subseteq$\\
   $\subseteq Im(A_1)+Im(A_2)+Im(A_3)\subseteq U_j$\\ (since $j\geqslant 4$), so \\$U_{j+1}=U_{j}+Im(A_{j+1})=U_{j}$\\
   
   Suppose now that $U_{k-1}=U_j.$ Then \\$Im(A_{k})=span \{v_{k,k-j}, \ v_{k,k-j+1}, \ v_{k,k-j+2}\}\subseteq $\\
   $\subseteq Im(A_{k-j})+Im(A_{k-j+1})+Im(A_{k-j+2})\subseteq U_{k-1}$\\
   (since $k-j\geqslant 1$ and $k-j+2\leqslant k-4+2=k-2$), so \\$U_{k}=U_{k-1}+Im(A_{k})=U_{k-1}=U_{j}$
   
  \end{proof}
  
    {\bf Theorem 5.2.13.} For  $n\geqslant 11,$ there are no {\it irreducible} matrix \\ representations  of $B_n$ of dimension $r\geqslant  n+1$ such that $rk(A_1)=3,$ 
$dim(Im(A_i)\cap Im(A_{i+1}))=0$  for all $i\in \mathbb{Z}_n,$ and \\
 $dim(Im(A_i)\cap Im(A_{j}))=1$ for  $||i-j||\geqslant 2,$ for all  $i, \ j\in\mathbb{Z}_n.$

 \begin{proof} Suppose that such irreducible representation $\rho$ exists.\\
 Consider subspaces $U_1, \ U_2, \ \dots , U_{n-1}$ defined by $U_1=Im(A_1)$ and $U_k=U_{k-1}+Im(A_k)$ for $k=2, \ 3, \ \dots , n-1.$\\ By Lemma 5.2.10 and Lemma 5.2.12, there exists a number $j,$\\ $  4\leqslant j \leqslant \frac{n}{2}+1,$ such that 
  $U_k=U_j$ for all $k=j+1, \ j+2, \dots , n-1.$\\
  By Lemma 5.2.11, $dim(U_j)\leqslant j+5\leqslant  \frac{n}{2}+1+5= \frac{n}{2}+6. $ \\
  Thus, $dim(U_k)=dim(U_j)\leqslant \frac{n}{2}+6$ for all $k=j+1, \ j+2, \dots , n-1.$\\ In particular, $dim(U_{n-1})\leqslant \frac{n}{2}+6.$\\
  
   On the other hand, 
   $U_{n-1}=Im(A_1)+Im(A_2)+\dots+Im(A_{n-1})=V,$ and 
   $n+1\leqslant  dim(V)= dim(U_{n-1})\leqslant \frac{n}{2}+6,$  which means $ \frac{n}{2}\leqslant 5, $ a contradiction with $n\geqslant 11.$
   
   \end{proof}
   
 Unfortunately, the method of the Theorem 5.2.13 does not work for $n=10$ in the case when $j=6.$ Thus, we have to consider this case separately.\\

  {\bf Theorem 5.2.14.}  For  $n=10,$ there are no {\it irreducible} matrix \\ representations  of $B_n$ of dimension $r\geqslant n+1$ such that $rk(A_1)=3,$ 
$dim(Im(A_i)\cap Im(A_{i+1}))=0$  for all $i\in \mathbb{Z}_n,$ and \\
 $dim(Im(A_i)\cap Im(A_{j}))=1$ for  $||i-j||\geqslant 2,$ $i, \ j\in\mathbb{Z}_n.$
  
\begin{proof} Suppose that such irreducible representation $\rho$ exists.\\
Consider subspaces $U_1, \ U_2, \ \dots , U_{9}$ defined by $U_1=Im(A_1)$  and\\ $U_k=U_{k-1}+Im(A_k)$ for $k=2, \ 3, \ \dots , 9.$\\ By Lemma 5.2.10 and Lemma 5.2.12, there exists a number $j,$\\$  4\leqslant j \leqslant \frac{n}{2}+1=6$ such that 
  $U_k=U_j$ for all $k=j+1, \ j+2, \dots , 9.$\\
  By Lemma 5.2.11, $dim(U_j)\leqslant j+5.$ If $j=4$ or $j=5,$ then \\
  $dim(U_j)\leqslant 10, $ and, in particular, $11\leqslant dim (V)=dim(U_{9})\leqslant \ 10,$ a contradiction.\\
  
  Suppose now that $j=6.$ In this case we have that for $j=4,\ 5,$\\ $X_{1,j}\neq Im(A_1).$ This means that $dim(X_{1,6})=3$ and \\$dim(X_{1,4})=dim(X_{1,5})=2.$ So, we have that \\
  $v_{1,5}\notin span\{v_{1,6},v_{1,7}\}=Y_{1,6}.$ \\

  By applying congugation by $T^2D$ to $X_{1,4}$ and $X_{1,5},$ we obtain that\\ $dim(X_{1,7})=dim(X_{1,8})=2.$ By Lemma 5.2.4, $dim(Y_{1,j})=2$ for $j\neq 0, \ 1, \ 2,$ and $j\neq 9.$ Thus, $X_{1,7}=X_{1,8}=Y_{1,6}=Y_{1,7}=Y_{1,8}.$ In particular,  $v_{1,6}\in span\{v_{1,8},v_{1,9}\}$ and $v_{1,5}\notin span\{v_{1,7},v_{1,8}\}.$\\
  In addition, by cyclic argument,   $v_{2,6}\notin span\{v_{2,8},v_{2,9}\}.$\\
  
  We claim that $span\{v_{1,6}\}$ is invariant under $B_{10}.$ \\
  First, since each of $A_1, \ A_3, \ A_4, \ A_6, \ A_8,$ and $A_9$ commute with both $A_1$ and $A_6,$ we have $A_iv_{1,6}\in span\{v_{1,6}\}$ for $i=1, \ 3, \ 4, \ 6, \ 8, \ 9.$\\
  
  Consider $A_2v_{1,6}.$ Since $A_2$ commutes with $A_6,$ we get $A_2v_{1,6}\in span \{v_{2,6}\}.$ On the other hand,  $v_{1,6}\in span\{v_{1,8},v_{1,9}\},$ and since $A_2$ commutes with both $A_8$ and $A_9,$ we have $A_2v_{1,6}\in span \{v_{2,8},v_{2,9}\}.$ Thus, \\ $A_2v_{1,6}\in span \{v_{2,6}\}\cap span \{v_{2,8},v_{2,9}\}=\{0\},$ since $v_{2,6}\notin span\{v_{2,8},v_{2,9}\},$ and hence, $A_2v_{1,6}=0.$\\
  
  By conjugating by $T^5,$ we get $A_7v_{6,1}=A_7v_{1,6}=0,$ and by conjugating by $T^2D$  we get $A_0v_{1,6}=0.$ Conjugating the last equation by $T^5$ gives $A_4v_{1,6}=0,$ which completes the proof that $span\{v_{1,6}\}$ is $B_{10}-$invariant. Thus, we get a contradiction with the irreducibility \\of $\rho.$
  
  \end{proof}

  \section{ The Main Theorem. }
  
  In this section we will summarize the main results of this paper.\\
  
  {\bf Theorem 6.1.} For $n\geqslant 10$ there are no irreducible complex \\representations of braid group $B_n$ on $n$ strings of dimension $n+1.$
  
 \begin{proof} By Theorem 5.1, every irreducible representation of dimension $n+1$ is equivalent to a tensor product of  a one-dimensional representation  and an irreducible representation of dimension $n+1$ and \\corank $3.$\\
  
  By Theorem 5.1.24, for every irreducible representation of dimension $n+1$ and corank $3,$ $Im(A_i)\cap Im(A_{i+1})=\{0\}$ for all $i\in \mathbb{Z}_n.$ Then by Lemma 5.2.1, $dim(Im(A_i)\cap Im(A_{j}))=1$ for all $||i-j||\geqslant 2, \  i, \  j \in \mathbb{Z}_n.$\\
  
  And finally,  by Theorems 5.2.13 and 5.2.14,  for $n\geqslant 10,$ there are no irreducible representations of dimension $n+1$ with $rk(A_1)=3,$   such that $dim(Im(A_i)\cap Im(A_{i+1}))=0$  for all $i\in \mathbb{Z}_n,$ and \\ 
  $dim(Im(A_i)\cap Im(A_{j}))=1$ for  $||i-j||\geqslant 2,$ $i, \ j\in\mathbb{Z}_n.$
  \end{proof}

  In addition,  without any extra effort, we get the following theorem.\\
  
   {\bf Theorem 6.2.} Let $\rho: B_n\to GL_{r}(\mathbb{C})$ be an {\it irreducible} matrix representation of $B_n$ for $n\geqslant 10.$ Then $rk(\rho(\sigma_i)-I)\neq 3, $ for all \\$i=1, \ 2, \dots, \ n-1.$
   
   \begin{proof} For $r\geqslant n+1,$ the statement follows from Theorems 5.1.24, 5.2.13 and 5.2.14.\\
   
   For $r\leqslant n,$ by \cite{Formanek} Theorem 22 and \cite{S} Theorem 6.1, $r=1$ or $n-2\leqslant r\leqslant n,$ and for $n-2\leqslant r\leqslant n,$ $\rho$ is equivalent to $\chi(y)\otimes\widehat{\rho},$ where  $\chi(y)$ is a one-dimensional representation, and $\widehat{\rho}$ is either a composition factor of a specialization of the reduced Burau representation, or a specialization of the standard representation.\\

    Clearly, if $r=1$ then $corank(\rho)\neq 3. $ For $n-2\leqslant r\leqslant n,$ if $y= 1,$ then $corank(\chi(y)\otimes\widehat{\rho}),$ is equal to 1 (when $r=n-2$ or $r=n-1$) or 2 (when $r=n$); and if $y\neq 1, $ then $corank(\chi(y)\otimes\widehat{\rho})\geqslant n-3\geqslant 7$  for $n\geqslant 10. $ 
    \end{proof}


\begin{thebibliography}{99}
  
  \bibitem{Burau} W. Burau, {\it Uber Zopfgruppen und gleichsinnig verdrillte Verkettungen},  Abh. Math. Sem. Univ. Hamburg {\bf 11} (1936), 179-186.

  \bibitem{Chow} W.-L. Chow, {\it On the algebraical braid group}, Ann. of Math. {\bf 49} (1948), 654-658.
  
  \bibitem{Formanek} E. Formanek, {\it Braid group representations of low degree}, Proc. London Math. Soc. {\bf 73} (1996), 279-322.
  
  \bibitem{FLSV} E. Formanek, W. Lee, I. Sysoeva, M. Vazirani, {\it The irreducible complex representations
of the braid group on $n$ strings of degree $\leqslant n$}, J. Algebra Appl. {\bf 2} (2003), 317–333.

\bibitem{Garside} F. A. Garside, {\it The Braid Group and Other Groups} Quart. J. Math. Oxford, {\bf 20} (1969), 235-254.

  \bibitem{Jones} V. F. R. Jones, {\it Braid groups, Hecke algebras and type ${\textrm{II}}_1$ factors}, from: “Geometric methods in operator algebras (Kyoto, 1983)”, Pitman Res. Notes Math. Ser. {\bf 123},
Longman Sci. Tech., Harlow (1986), 242–273.

\bibitem{Jones87}  V. F. R. Jones, {\it Hecke Algebra Representations of Braid Groups and Link Polynomials}, Ann. of Math. {\bf 126} (1987), 335-388.

\bibitem{LarsenRowell} M. Larsen, E. Rowell {\it  Unitary braid representations with finite image}, Algebr. Geom. Topol. {\bf 8} (2008), No. 4, 2063--2079.

\bibitem{Lawrence} R. Lawrence, {\it Homological representations of the
Hecke algebra}, Comm. Math. Phys. {\bf 135}, No. 1, (1990), 141-191. 

\bibitem {Lee} W. Lee, {\it Representations of the braid group $B_4$} , J. Korean Math. Soc. {\bf 34}, No. 3(1997), 673-693.

  \bibitem{S} I. Sysoeva, {\it Dimension $n$ representations of the braid group on $n$ strings}, J. Algebra {\bf 243}, (2001) 518–538.
  
 
 \bibitem{TongYangMa} Dian-Min Tong, Shan-De Yang, and Zhong-Qi Ma, {\it A new class of representations of
braid groups}, Comm. Theoret. Phys. {\bf 26}, No. 4 (1996), 483-486. 

\bibitem{Zinno}  M. Zinno, {\it On Krammer's representation of the braid group}, Math. Ann. {\bf 321}, No. 1 (2001), 192-211. 

  
  \end{thebibliography}
\end{document}